\documentclass[10pt]{article}
\usepackage{palatino}
\usepackage{amsmath}
\usepackage{amssymb}
\usepackage{amscd}
\usepackage{latexsym}
\usepackage{theorem}
\usepackage{doublespace}
\usepackage{epsfig}
\usepackage{psfrag}
\usepackage{amsfonts}
\usepackage{enumerate}
\usepackage{euscript}
\usepackage{array}
\usepackage{bbm}

\newtheorem{theorem}{Theorem}[section]
\newtheorem{lemma}[theorem]{Lemma}
\newtheorem{proposition}[theorem]{Proposition}
\newtheorem{corollary}[theorem]{Corollary}

{\theorembodyfont{\rmfamily}
\theoremstyle{plain}
\newtheorem{definition}[theorem]{Definition}

\newtheorem{remark}[theorem]{Remark}

}

\numberwithin{equation}{section}

\input xy
 \xyoption{all}

\newcommand{\hs}{{\hspace{3mm}}}
\newcommand{\hsm}{{\hspace{1mm}}}

\newcommand{\gl}{{\mathfrak gl}}

\renewcommand{\H}{{\mathbb{H}}}
\newcommand{\C}{{\mathbb{C}}}

\newcommand{\Z}{{\mathbb{Z}}}

\newcommand{\R}{{\mathbb{R}}}

\newcommand{\into}{{\hookrightarrow}}

\newcommand{\qed}{\hfill \mbox{$\Box$}\medskip\newline}
\newenvironment{proof}{\noindent {\bf Proof:}}{\qed \par}

\newenvironment{question}{\noindent {\bf Question:}}{\newline}

\newcommand{\Id}{\mathbbm{1}}
\newcommand{\Hn}{{\mathbb H}^n}
\newcommand{\g}{{\mathfrak g}}

\renewcommand{\sp}{{\mathfrak sp}}
\renewcommand{\t}{{\mathfrak t}}
\renewcommand{\u}{{\mathfrak u}}

\newcommand{\unH}{{\mathfrak u}(n,\H)}
\newcommand{\unHs}{{\mathfrak u}(n-1,\H)}
\newcommand{\unC}{{\mathfrak u}(n,\C)}

\newcommand{\unCs}{{\mathfrak u}(n-1,\C)}
\newcommand{\tdZ}{{\mathfrak t}^{*}_{\Z}}
\newcommand{\intwc}{(\t^{*}_{+})_{0}}

\renewcommand{\P}{{\mathbb P}}

\newcommand{\Ol}{{\mathcal O}_{\lambda}}
\newcommand{\Ored}{{\mathcal O}_{\lambda} \hsm /\!/_{\!\mu} \hsm U(n-1,\H)}
\newcommand{\mult}{M_{\lambda}^{\mu}}
\newcommand{\mmod}{/\!/}

\newcommand{\Ucent}{U(\sp(2n,\C))^{\sp(2(n-1),\C)}}

\newcommand{\sptyt}{{\mathcal G}_{2}/{\mathcal H}_2}

\newcommand{\kerPhi}{{\rm ker}(d\overline{\Phi}_{n-1})}

\begin{document}
\begin{spacing}{1.1}

\noindent
{\LARGE \bf The symplectic geometry of the Gel'fand-Cetlin-Molev basis
for representations of $Sp(2n,\C)$}
\bigskip\\
{\bf Megumi Harada}\footnote{{\tt megumi@math.toronto.edu}.
\newline \mbox{~~~~} {\it MSC 2000 Subject Classification}: 
Primary: 70H06 \hspace{0.1in} Secondary: 37J15, 17B10
\newline \mbox{~~~~} {\it Keywords:} 
integrable systems, coadjoint orbits, symplectic group}  \\ 
Department of Mathematics, University of Toronto, 
Toronto, Ontario, Canada M5S 3G3 \smallskip \\
\bigskip

\bigskip

{\small
\begin{quote}
\noindent {\em Abstract.}  Gel$'$fand and Cetlin constructed in the
1950s a canonical basis for a finite-dimensional representation
$V(\lambda)$ of $U(n,\C)$ by successive decompositions of the
representation by a chain of subgroups \cite{GelTse2, GelTse1}.  
Guillemin and Sternberg
constructed in the 1980s the {\em Gel$'$fand-Cetlin integrable system}
on the coadjoint orbits of $U(n,\C)$, which is the symplectic
geometric version, via geometric quantization, of the
Gel$'$fand-Cetlin construction. (Much the same construction works for
representations of $SO(n,\R)$.)  A. Molev \cite{M4} in 1999
found a Gel$'$fand-Cetlin-type basis for representations of the
symplectic group, using essentially new ideas. An important new role
is played by the Yangian $Y(2)$, an infinite-dimensional Hopf algebra,
and a subalgebra of $Y(2)$ called the twisted Yangian $Y^{-}(2)$. In
this paper we use deformation theory to give the analogous
symplectic-geometric results for the case of $U(n,\H)$, i.e. we
construct a completely integrable system on the coadjoint orbits of
$U(n,\H)$. We call this the {\em Gel$'$fand-Cetlin-Molev integrable
system.}
\end{quote}
}
\bigskip

\tableofcontents

\section{Introduction}\label{sec:intro}

Symplectic geometry and representation theory can related by the
theory of {\em geometric quantization}. In this theory, a symplectic
manifold $M$ equipped with a Hamiltonian $G$-action has an associated
linear representation $V$ of $G$, and symplectic reductions by $G$ of
$M$ translate to taking $G$-isotypic components $V^{\lambda}$ of
$V$. This correspondence has served as an underlying theme in much
work in modern symplectic geometry, and the present paper is no
exception. 

There are two parallel theories which, as a pair, serve as motivation
for the work in this paper. These are the {\em Gel'fand-Cetlin basis}
for representations of $U(n,\C)$, and the {\em Gel'fand-Cetlin
integrable system} on coadjoint orbits of $U(n,\C)$. The names are no
coincidence: via geometric quantization, the Gel'fand-Cetlin system on
a coadjoint orbit can be seen to be the symplectic-geometric analogue
of the Gel'fand-Cetlin basis for an appropriate representation. This
parallel works out beautifully in the case of $U(n,\C)$ or $SO(n,\R)$,
but not for other groups. (For simplicity in the discussion below, we
refer only to the group $U(n,\C)$.)

In both of these parallel theories, the underlying goal is to produce
a ``large'' torus action on the relevant space (in the case of
representation theory, a vector space, and in the case of symplectic
geometry, a symplectic manifold). 
On the representation-theoretic side, the torus 
is ``maximal'' in the sense that it completely decomposes the
representation into 1-dimensional
eigenspaces. On the symplectic side, the torus is ``maximal'' in the
well-known sense that the torus is half the dimension of the
symplectic manifold, thus making an integrable system. 

A finite-dimensional irreducible representation $V_{\lambda}$ of
$U(n,\C)$, when considered as a representation of the subgroup
$U(n-1,\C)$, decomposes with multiplicity 1.  The representation
$V_{\lambda}$ can be decomposed successively by a chain of subgroups
\[U(1,\C) \subset U(2,\C) \subset \ldots \subset U(n-1,\C).\] 
Since the final subgroup \(U(1,\C)\) is abelian, and because
of the multiplicity-free decomposition at each step, one obtains
a canonical (up to a choice of this chain of subgroups) basis 
for any finite-dimensional $U(n,\C)$-representation $V_{\lambda}$. 
This basis is called the {\em Gel'fand-Cetlin basis} for
$V_{\lambda}$ \cite{GelTse1}. Its construction is briefly summarized in
Section~\ref{subsec:GCforUnC}.

Guillemin and Sternberg showed in the 1980s \cite{GS83} that a special
set of functions on a coadjoint orbit $\Ol$ of the unitary group
$U(n,\C)$ give the maximal possible number of Poisson-commuting
functions on $\Ol$. In other words, they show that the coadjoint orbit
$\Ol$ is a completely integrable system.  (The coadjoint
orbits are {\em not} toric varieties; this is because the
aforementioned functions are smooth only on an open dense subset of
$\Ol$.) The coadjoint orbit, equipped with these Poisson-commuting
functions, is called the {\em Gel'fand-Cetlin system} on $\Ol$. By
geometric quantization, the existence of this maximal set of
Poisson-commuting functions is the geometric analogue of the multiplicity-free
decomposition of $V_{\lambda}$. This integrable system is explained in
Section~\ref{subsec:GCSystemforUnC}. 

For other groups, finding a Gel$'$fand-Cetlin basis proves to be much
more difficult. In particular, for the case of $U(n,\H)$ (the compact
form of $Sp(2n,\C)$), a
difficulty arises in that the finite-dimensional irreducible
representations $V_{\lambda}$ of $U(n,\H)$ decompose {\em with
multiplicity} as representations of $U(n-1,\H)$. Similarly, from the
symplectic-geometric standpoint, the symplectic reductions of $\Ol$ by
$U(n-1,\H)$ are not just points (as they are for the $U(n,\C)$ case),
but are nontrivial symplectic manifolds. 

Nevertheless, A. Molev
\cite{M4} found a Gel$'$fand-Cetlin-type basis for
finite-dimensional irreducible representations of $U(n,\H)$, which are
constructed in the spirit of the original work of
Gel$'$fand-Cetlin. His methods required the use of new tools,
including an infinite-dimensional algebra called the Yangian. Molev's
theorems are recounted in Section~\ref{subsec:GCforUnH}. 

This history is summarized in
the table below. The essence of this paper is to answer the following

\begin{center}
\parbox{5in}{
\begin{question}
What is the ``Gel$'$fand-Cetlin-type'' integrable system on coadjoint
orbits of $U(n,\H)$ corresponding (via geometric quantization) to
Molev's canonical bases of the representations? In other words, what
goes in the bottom right-hand corner of the table below? 
\end{question}}
\end{center}

\bigskip
\bigskip

\begin{tabular}{|m{1.5cm}|m{6cm}|m{7cm}|} 
\hline 
& & \\
& {\bf SYMPLECTIC GEOMETRY} & {\bf REPRESENTATION THEORY} \\ 
& & \\
\hline 
& & \\
$U(n,\C)$ & Gel$'$fand-Cetlin canonical basis for
finite-dimensional irreducible representations $V(\lambda)$  &
Gel$'$fand-Cetlin integrable system on coadjoint orbits $\Ol$ \\ 
 & \hspace{1cm} (Gel$'$fand-Cetlin, 1950) & \hspace{1cm} (Guillemin-Sternberg, 1983) \\ 
& & \\
\hline 
& & \\
$U(n,\H)$ & Gel$'$fand-Cetlin-type canonical basis for finite-dimensional
irreducible representations $V(\lambda)$ & \hspace{3cm} {\Huge\bf ???}  \\ 
 & \hspace{1cm} (Molev, 1999) &   \\ 
& & \\
\hline 
\end{tabular}

\bigskip
\bigskip

It turns out that the answer to this question, which is the analogous theorem to
Guillemin and Sternberg's for the case of the group \(G = U(n,\H),\)
is remarkably easy to state. 
For simplicity, we always assume that
$\Ol$ is a {\em generic} coadjoint orbit of $U(n,\H)$. 
For convenience, we first state the result in terms of an
intermediate geometric object, namely the
symplectic reductions of $\Ol$ by the subgroup $U(n-1,\H)$.

\begin{theorem}\label{thm:Oredfunctions}
Let \({\mathcal O}_{\lambda} \cong U(n,\H)/T^n\) be a coadjoint orbit of 
$U(n,\H)$. Let $\Psi$ be the $n$-th component of the $T^n$ moment map
on $\Ol$. Let $g_{n,m},$ for $1 \leq {m} \leq n-1$ be defined by 
\[
g_{n,m}([A]) = |a_{n{m}}|^2, 
\]
where \(A = (a_{ij}) \in U(n,\H)\) and 
\(a_{ij} \in \H.\) Then the functions $\{g_{n,m}\}_{m=1}^{n-1}$ and $\Psi$ 
descend to a completely integrable system on 
the reduced space \({\mathcal O}_{\lambda} \mmod\!_{\mu} U(n-1,\H)\)
for $\mu$ a regular value.  
\end{theorem} 

Using these functions, plus an inductive construction,
one can show that the original coadjoint orbit 
$\Ol$ is also an integrable system. 

\begin{theorem}\label{thm:int_sys}
Let $\Ol$ be a (generic) coadjoint orbit of $U(n,\H)$. 
Then there exists a maximal set of Poisson-commuting
functions on $\Ol$, making it an integrable system. 
\end{theorem}

We call this the Gel$'$fand-Cetlin-Molev integrable system. The proofs
of these main results are in Section~\ref{sec:GCM}.
Although the formul{\ae} for the $g_{n,m}$ given
above are remarkably simple, it is nevertheless instructive to reveal
how these formul{\ae} were derived, via the theory of deformation
quantization, from Molev's work. In particular, the 
technical heart of the derivation of the formul{\ae} used in
Theorem~\ref{thm:Oredfunctions} lies in the remarkable fact that
certain algebra automorphisms of the Yangian $Y(2n)$ degenerate,
in the classical limit, to be {\em trivial}. This story is presented in
Section~\ref{sec:limits}.

\bigskip

\noindent{\bf Acknowledgements.} 
The author thanks Allen Knutson for suggesting this problem, as
well as
for help and encouragement throughout. Many thanks also to Alexandre Molev for
answering questions about twisted Yangians.

\section{History}\label{sec:history}

\subsection{The Gel$'$fand-Cetlin basis for $U(n,\C)$-representations}\label{subsec:GCforUnC}

We now briefly recall the construction of the Gel$'$fand-Cetlin basis \cite{GelTse1}.
Let $U(n,\C)$ be the standard unitary group, 
acting on $\C^n$ equipped with the standard hermitian form 
and standard orthonormal basis denoted by \(\{e_1, \ldots, e_n\}.\) 

We consider the chain of subgroups 
\begin{equation}\label{eq:chain}
U(1,\C) \subset U(2,\C) \subset \ldots \subset U(n-1,\C) 
\subset U(n,\C),
\end{equation}
where $U(k,\C)$ is the subgroup of $U(n,\C)$
fixing \(\{e_{k+1}, \ldots, e_n\}.\) 

Let \(V(\lambda)\) be the irreducible representation of $U(n,\C)$ of
highest weight \(\lambda = (\lambda_1 \geq \lambda_2
\geq \ldots \geq \lambda_n) \in (\t^n)^{*}_{\Z}.\) 
Considered as a $U(n-1,\C)$-representation, $V(\lambda)$ may
not be irreducible.  Indeed, it will decompose as
\[
V(\lambda) \cong \bigoplus_{\mu} V(\lambda)^{\mu} \quad \mbox{ as }
U(n-1,\C) \mbox{ representations, } 
\]
where \(W^{\mu}\) denotes the $\mu$-isotypic component. 
We denote by \(\mu=(\mu_1,
\ldots, \mu_{n-1}) \in (\t^{n-1})^{*}_{\Z}\)
a dominant weight for $U(n-1,\C)$. 

The $\mu$-isotypic component of $V(\lambda)$ may, a priori, contain many 
copies of the irreducible representation $V(\mu)$ of $U(n-1,\C)$. 
In other words, we may also write 
\[
V(\lambda)^{\mu} \cong M_{\lambda}^{\mu} \otimes V(\mu),
\hs \quad \mbox{ as }
U(n-1,\C) \mbox{ representations,}
\]
where $U(n-1,\C)$ acts trivially on $M_{\lambda}^{\mu}$, 
and the multiplicity space $M_{\lambda}^{\mu}$ is 
the subspace of high-weight vectors in the $\mu$-isotypic component, i.e.
\[
M_{\lambda}^{\mu} := (V(\lambda)^{\mu})^+.
\]
Thus we have 
\[
V(\lambda) \cong \bigoplus_{\mu} (M_{\lambda}^{\mu} 
\otimes V(\mu)) \quad \mbox{ as }
U(n-1,\C) \mbox{ representations.}
\]
The following two facts are crucial. 
First, it turns out that 
\(\text{dim}(M_{\lambda}^{\mu}) \not = 0\)
if and only if 
\begin{equation}\label{eq:UnCinequalities}
\lambda_1 \geq \mu_1 \geq \lambda_2 \geq
\ldots \lambda_{n-1} \geq \mu_{n-1} \geq \lambda_n.
\end{equation} 
Second, for $\mu$ that do appear in the decomposition, \(\text{dim}(M_{\lambda}^{\mu})=1.\) 
In other words, the decomposition
is multiplicity-free. 

\begin{figure}[h]
\begin{center}
\psfrag{VL}{$V(\lambda)$}
\psfrag{VLm1}{$V(\lambda)^{\mu_1}$}
\psfrag{VLm2}{$V(\lambda)^{\mu_2}$}
\psfrag{VLm3}{$V(\lambda)^{\mu_3}$}
\psfrag{VLm1n1}{$(V(\lambda)^{\mu_1})^{\nu_1}$}
\psfrag{VLm1n2}{$(V(\lambda)^{\mu_1})^{\nu_2}$}
\psfrag{VLm1n3}{$(V(\lambda)^{\mu_1})^{\nu_3}$}
\psfrag{C}{${\mathbb C}$}
\psfrag{Tn-1}{$T^{n-1}$}
\psfrag{Tn-2}{$T^{n-2}$}
\psfrag{T1}{$T^1$}
\epsfig{figure=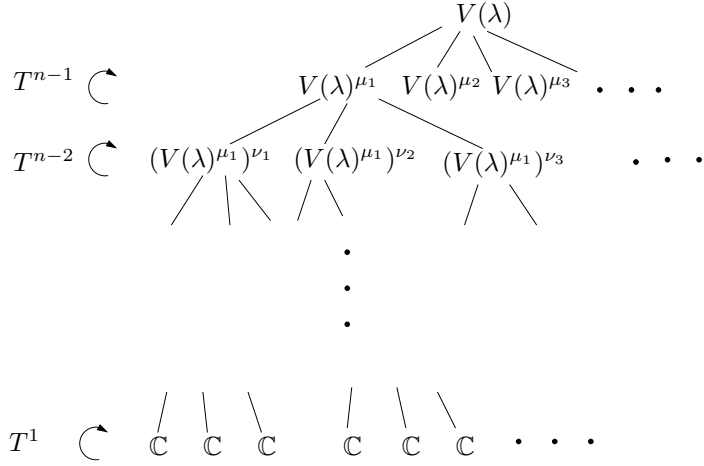, height=6cm}
\end{center}
\begin{center}
\parbox{5in}{
\caption{At each step in the successive decomposition, we get a
  $T^{n-k}$-action for appropriate $k$. Since at each step
  the $T^{n-k}$ act as scalars on each isotypic component, all the
  $T^{n-k}$ actions commute with each other. Hence we get a
  $T^{n(n-1)/2}$-action.}\label{fig:Vdecomp}
}
\end{center}
\end{figure}

Recall that the maximal torus $T^{n-1}$ of $U(n-1,\C)$,
the subgroup of diagonal matrices in $U(n-1,\C)$, 
naturally acts on $V(\lambda)$ by restriction. 
We now define a new $T^{n-1}$-action on $V(\lambda)$, 
different from, though related to, the above action. 
Namely, we define the new torus $T^{n-1}$ to act 
on each $\mu$-isotypic component 
as scalar matrices, so that each non-zero vector
behaves as a $T^{n-1}$-weight vector of weight $\mu$.
In other words, for \(v \in V(\lambda)^{\mu},\) we define
\[
t \cdot v := \mu(t)v.
\]
We call this the Gel$'$fand-Cetlin $T^{n-1}$-action on
$V(\lambda)^{\mu}$.  We now repeat this process, using the subgroup
$U(n-2,\C)$ of $U(n-1,\C)$ in~\eqref{eq:chain}.
The key observation now is that the new $T^{n-2}$-action, defined in the
same fashion, commutes with the $T^{n-1}$-action defined
previously.  This is because the $T^{n-1}$ acts by scalar matrices on
each component $V(\mu)$. Since $V(\lambda)$ is a sum of the $V(\mu)$,
the two tori commute on $V(\lambda)$. Thus we have now a $T^{n-1}
\times T^{n-2}$-action on $V(\lambda)$.  By continuing this process
for the whole chain of subgroups, we obtain at each step a
$T^{k}$-action commuting with the previous $T^{k+1}$. Hence when we 
reach the last subgroup $U(1,\C)$, we have obtained a \(T^{n-1} \times T^{n-2}
\times \cdots \times T^2 \times T^1 \cong T^{n(n-1)/2}\) action on
$V(\lambda)$. This decomposition is schematically illustrated in
Figure~\ref{fig:Vdecomp}.

Now consider the decomposition of $V(\lambda)$ into $T^{n(n-1)/2}$-weight
spaces. Since the decomposition is multiplicity-free
at each step, and since the last group $U(1,\C) \cong S^1$ is 
abelian, each $T^{n(n-1)/2}$-weight space is one-dimensional. Hence
this torus completely decomposes the representation $V(\lambda)$
into 1-dimensional subspaces, thus providing (up to a choice of 
scalar in each weight space) a canonical basis for $V(\lambda)$. 
This basis is called the {\em Gel$'$fand-Cetlin basis} for $V(\lambda)$.

\begin{figure}[h]
\begin{center}
\psfrag{L1}{$\lambda_1$}
\psfrag{L2}{$\lambda_2$}
\psfrag{L3}{$\lambda_3$}
\psfrag{L4}{$\lambda_4$}
\psfrag{Ln-2}{$\lambda_{n-2}$}
\psfrag{Ln-1}{$\lambda_{n-1}$}
\psfrag{Ln}{$\lambda_n$}
\psfrag{M1}{$\mu_1$}
\psfrag{M2}{$\mu_2$}
\psfrag{M3}{$\mu_3$}
\psfrag{Mn-2}{$\mu_{n-2}$}
\psfrag{Mn-1}{$\mu_{n-1}$}
\psfrag{N1}{$\nu_1$}
\psfrag{N2}{$\nu_2$}
\psfrag{Nn-2}{$\nu_{n-2}$}
\psfrag{R1}{$\rho_1$}
\epsfig{figure=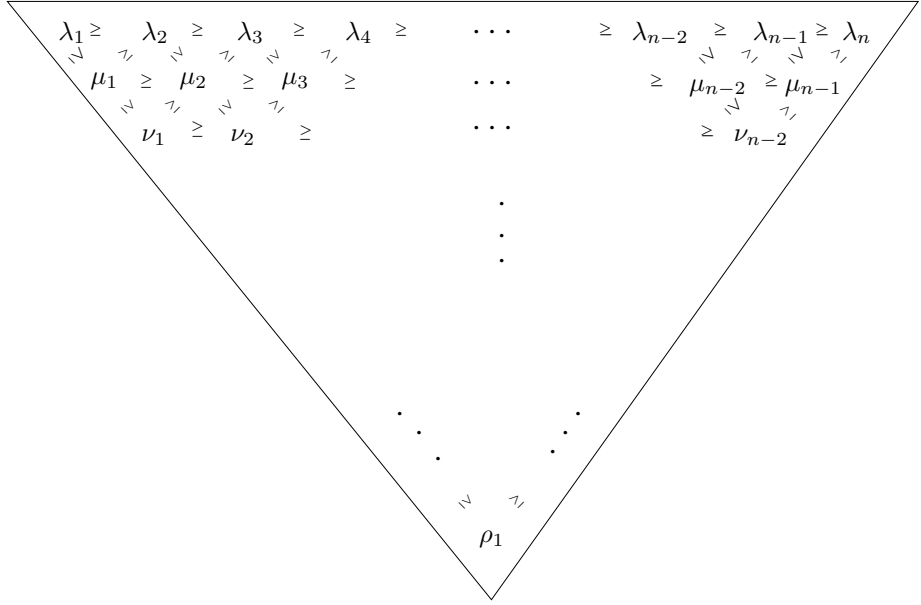, height=8cm}
\end{center}
\begin{center}
\parbox{5in}{
\caption{The integer arrays parametrizing the Gel$'$fand-Cetlin basis
for $U(n,\C)$ representations. The top row is fixed, given by the highest weight of the
irreducible representation $V(\lambda)$.} \label{fig:GCforUnC}
}
\end{center}
\end{figure}

  It is pleasant to note that this construction also gives a {\em
  combinatorial} algorithm for counting the dimension of any
  finite-dimensional irreducible represenation of $V(\lambda)$. As a
  result of the facts that \(\rm{dim}(M_{\lambda}^{\mu}) \neq 0\)
  if and only if the inequalities~\eqref{eq:UnCinequalities} are
  satisfied and that the decomposition is multiplicity-free, the
  dimension of $V(\lambda)$ is given by the number of integer fillings of the triangle
  in Figure~\ref{fig:GCforUnC} with specified top row, satisfying the
  given inequalities.

\subsection{The Gel$'$fand-Cetlin integrable system for $U(n,\C)$}\label{subsec:GCSystemforUnC}

We now explain the symplectic-geometric side of the Gel$'$fand-Cetlin
story for $U(n,\C)$ \cite{GS83}. 
By geometric quantization, the symplectic-geometric object
corresponding to an irreducible representation $V(\lambda)$ is the
coadjoint orbit $\Ol$ of $U(n,\C)$ through the point \(\lambda \in
\tdZ \subseteq \t^{*} \subseteq \g^{*}\) \cite{Kos65}. Here we use the Killing form
to identify \(\g \cong \g^{*},\) and will think of $\t^{*}$ as a
subspace of $\g^{*}$. 
Such a coadjoint orbit $\Ol$ is a symplectic
manifold, generically of (real) dimension $n(n-1)$.
Hence, the result analogous to the
existence of a Gel$'$fand-Cetlin basis for $V(\lambda)$ will be the
existence of a $\frac{n(n-1)}{2}$-dimensional torus
acting in a Hamiltonian fashion on $\Ol$.

We construct this large torus action on $\Ol$ as follows. As in the
representation-theoretic construction above, we will consider the
actions of the subgroups $U(n-k,\C)$ in the chain~\eqref{eq:chain},
and use the tori $T^{n-k}$ in each $U(n-k,\C)$.  The coadjoint orbit
$\Ol$ is a Hamiltonian $U(n,\C)$-manifold with moment map \(\Phi: \Ol
\into \g^{*}\) the inclusion.  We may restrict to the subgroups
$U(n-k,\C)$, which also act in a Hamiltonian fashion on $\Ol$ with
moment maps given by composing $\Phi$ with the projections
\[\pi_{n-k}: \unC^{*} \to \u(n-k,\C)^{*}.\] 
Thus we have a collection of moment maps \(\Phi_{n-k} := \pi_{n-k}
\circ \Phi,\) as shown below:  
\begin{equation}\label{eq:Phik}
\xymatrix{
\Ol \ar @{^{(}->}[r]^(0.4){\Phi} \ar @/^2pc/[rr]_{\Phi_{n-1}} \ar @/^3pc/[rrr]^{\Phi_{n-2}} 
& \unC^{*} \ar[r] & 
 \unCs^{*} \ar[r] &  \u(n-2,\C)^{*} \ar[r] &  \cdots  
}
\end{equation}
given by successively projecting onto smaller \(\u(n-k,\C)^{*}.\) 

The large torus action is obtained by using the {\em Thimm trick} on
each $U(n-k,\C)$ moment map. By taking the projections
to the positive Weyl chamber for each moment map $\Phi_k$, 
we get a sequence of maps to smaller and smaller Weyl chambers.
\begin{equation}\label{eq:ThimmforUnC}
\xymatrix{
\Ol  \ar @{^{(}->}[r] & \unC^{*} \ar[r]  & 
\unCs^{*} \ar[r] \ar[d]^{/U(n-1,\C)} & \u(n-2,\C)^{*} \ar[r] 
\ar[d]^{/U(n-2,\C)} & \cdots \\
   &   & (\t^{n-1})_{+}^{*}   & (\t^{n-2})^{*}_{+}  & \cdots 
}
\end{equation}
In linear-algebraic terms,
each of these projections is given by diagonalizing a matrix in
$\u(k,\C)^* \cong \u(k,\C)$ and reading off the diagonal entries
(arranged to be in non-increasing order). The first
projection from \(\u(n,\C)^* \to (\t^n)^*_{+}\) is omitted since it is
trivial when restricted to the fixed coadjoint orbit $\Ol$. 

By the Thimm trick, these
functions to the positive Weyl chambers give rise to a torus action
on an open dense set in $\Ol$. The action of the Thimm torus on $M$
is, heuristically, given by ``moving to the symplectic slice, acting
by the (usual) torus, then moving back.''  This is analogous to having
the tori $T^{n-k}$ act as {\em scalars} on the whole $V(\mu)$ instead
of just on the high-weight vectors. We will call these the {\em
Gel$'$fand-Cetlin tori $T^{n-k}$} acting on $\Ol$ for each $k$, in
analogy with the representation-theoretic situation.  These tori also
commute with one another since their moment maps are Casimirs
(i.e. constant on symplectic leaves).  Thus, there is an action of
\(T^{n-1} \times T^{n-2} \times \cdots \times T^1 \cong
T^{n(n-1)/2},\) as advertised, on $\Ol$.  Since this is the maximal
possible dimension of a torus acting Hamiltonianly on $\Ol$, this is
called the {\em Gel$'$fand-Cetlin integrable system}.
Note that the $T$-action on $\Ol$ coming from the maximal torus
\(T \subset U(n,\C)\) is a sub-torus action of the Gel$'$fand-Cetlin
torus action. The analogous statement will {\em not} be true in the $U(n,\H)$ case.

We summarize, in the ``Rosetta Stone'' below, the correspondences
between specific objects arising in the two related
constructions. 
(Here ``G-C'' stands for ``Gel'fand-Cetlin.'')

\bigskip

\begin{tabular}{|m{7cm}|m{7cm}|}
\hline
\hline
{\bf Symplectic Geometry } & {\bf Representation Theory} \\
\hline
\hline 
coadjoint orbit $\Ol$ &  irreducible representation $V(\lambda)$ \\
\hline
\(U(n,\C)\) action on $\Ol$  & $U(n,\C)$  action on $V(\lambda)$ \\
\hline
\(U(n-1,\C) \subset U(n,\C)\) action on $\Ol$ &  \(U(n-1,\C)
\subset U(n,\C)\)  action on $V(\lambda)$ \\
\hline
\(\Phi^{-1}({\mathcal O}_{\mu})/\text{Stab}(\mu)\) &  $\mu$-isotypic component 
$V(\lambda)^{\mu}$ \\
\hline 
symplectic slice \(S = \Phi^{-1}(\intwc)\) &  high-weight
vectors \(V(\lambda)^{+}\) \\
\hline
Thimm torus action \(T^{n-1}\)  on $\Ol$ & G-C torus action 
\(T^{n-1}\) on $V(\lambda)$ \\
\hline 
symplectic reduction \({\mathcal O}_{\lambda} \hsm /\!/_{\!\mu} \hsm U(n-1,\C)\) & multiplicity space 
\(((V(\lambda))^{\mu})^{+} \cong M_{\lambda}^{\mu}\) \\
\hline
\hline
\end{tabular}

\bigskip
\noindent The last correspondence between the symplectic reduction and the
multiplicity space is the content of the
``quantization-commutes-with-reduction'' theorem in \cite{GS82}. In
particular, in the case of the Gel$'$fand-Cetlin system for $U(n,\C)$,
the fact that the multiplicity spaces $M_{\lambda}^{\mu}$ are
dimension 1 correspond to the symplectic geometric fact that the
symplectic reductions $\Ored$ are just {\em points.}

\subsection{The Gel$'$-fand-Cetlin basis for $U(n,\H)$-representations}\label{subsec:GCforUnH}

Molev's construction in \cite{M4} of the analogous Gel$'$fand-Cetlin
basis for finite-dimensional irreducible representations $V(\lambda)$
is phrased in terms of the complex group $Sp(2n,\C)$, and in this
section we do the same\footnote{We will abuse notation throughout and
use the same notation for objects associated to $U(n,\C)$ and the
corresponding objects associated to $U(n,\H)$. We hope the context
will make clear the group under discussion.}.  We will now briefly
recount his results, following his convention of using the complex
group $Sp(2n,\C)$.

As in the construction of the Gel$'$fand-Cetlin basis for
$GL(n,\C)$ representations, we first fix a choice of chain of
subgroups
\begin{equation}\label{spchain}
Sp(2,\C) \subset Sp(4,\C) \subset \cdots \subset Sp(2(n-1),\C) 
\subset Sp(2n,\C).
\end{equation}
Let $V(\lambda)$ be a finite-dimensional irreducible representation 
of $Sp(2n,\C)$, of highest weight $\lambda \in \tdZ$. Since the Weyl
group of $Sp(2n,\C)$ is the group of {\em signed} permutations, we follow
Molev's conventions in \cite{M4} and choose the 
positive Weyl chamber so that \(\lambda
= (\lambda_1, \lambda_2, \ldots, \lambda_n) \in \Z^n,\) and 

\begin{equation}\label{sp_wcineq}
0  \geq  \lambda_1 \geq \lambda_2 \geq \cdots \geq \lambda_{n-1} \geq \lambda_n.
\end{equation}
We now restrict to the action of the subgroup $Sp(2(n-1),\C)$ on
$V(\lambda)$. As in the $GL(n,\C)$ case, we obtain a decomposition
\begin{equation}\label{spdec}
V(\lambda) \cong \bigoplus_{\mu}
(M_{\lambda}^{\mu} \otimes V(\mu)) \quad \mbox{ as }
Sp(2(n-1),\C) \mbox{ representations, } 
\end{equation}
where $V(\mu)$ is a $Sp(2(n-1),\C)$-irreducible representation of 
highest weight $\mu$. 

The main difficulty in the $Sp(2n,\C)$ case is that the decomposition
above in~\eqref{spdec} is {\em not} multiplicity-free,
i.e. \(\text{dim}(M_{\lambda}^{\mu})\) is not necessarily equal to
1. Thus, the Thimm torus $T^{n(n-1)/2}$, constructed exactly as in the
case of $U(n,\C)$, acting on the decomposition by the
chain~\eqref{spchain} will decompose $V(\lambda)$ into smaller
subspaces, but {\em not} into 1-dimensional pieces. Hence, in order to
obtain a complete decomposition, we must find an additional action on
the multiplicity spaces $M_{\lambda}^{\mu}$.

We need not look very far to find an algebra acting on
$M_{\lambda}^{\mu}$. Since the multiplicity spaces $\mult$ are the
subspaces of high-weight vectors $(V(\lambda)^{\mu})^{+}$, the
centralizer \(\Ucent\) of \(\sp(2(n-1),\C)\) in \(U(\sp(2n,\C)),\)
acts on each of the highest-weight spaces \(V(\lambda)_{\mu}^{+}.\) In
fact, it is known that it acts {\em irreducibly} \cite[Section
9.1]{Dix74}. However, it is difficult to find explicitly the weights
of vectors in $\mult$ for an appropriate commuting subalgebra of
$\Ucent$, thus extending the Thimm torus action.

\begin{figure}[h]
\begin{center}
\psfrag{0}{$0$}
\psfrag{L1}{$\lambda_1$}
\psfrag{L2}{$\lambda_2$}
\psfrag{L3}{$\lambda_3$}
\psfrag{Ln}{$\lambda_n$}
\psfrag{L1p}{$\lambda'_1$}
\psfrag{L2p}{$\lambda'_2$}
\psfrag{L3p}{$\lambda'_3$}
\psfrag{L4p}{$\lambda'_4$}
\psfrag{Lnp}{$\lambda'_n$}
\psfrag{M1}{$\mu_1$}
\psfrag{M2}{$\mu_2$}
\psfrag{M3}{$\mu_3$}
\psfrag{Mn-1}{$\mu_{n-1}$}
\psfrag{M1p}{$\mu'_1$}
\psfrag{M2p}{$\mu'_2$}
\psfrag{M3p}{$\mu'_3$}
\psfrag{Mn-1p}{$\mu'_{n-1}$}
\psfrag{R1}{$\rho_1$}
\psfrag{R1p}{$\rho'_1$}
\epsfig{figure=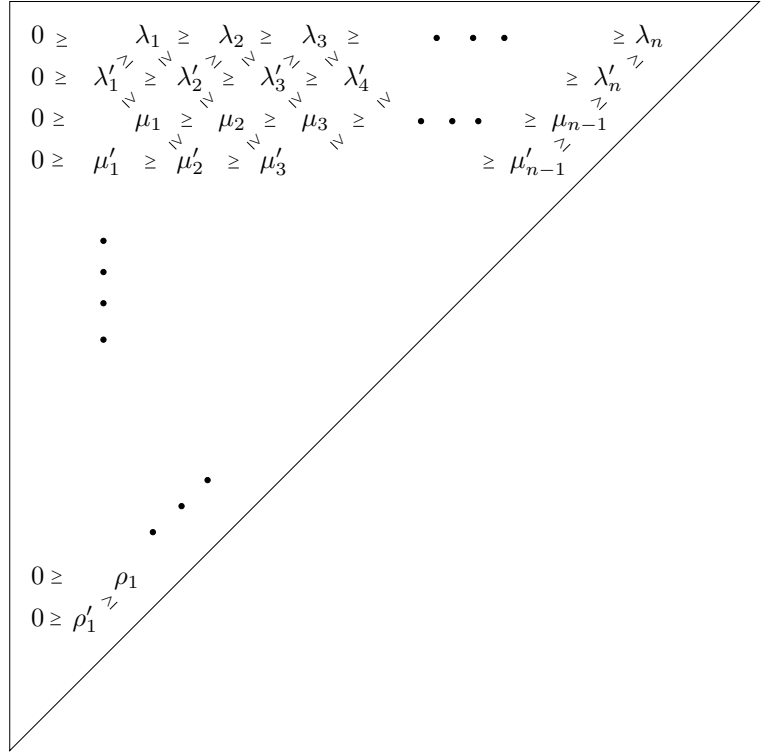, height=10cm}
\end{center}
\begin{center}
\parbox{5in}{
\caption{The integer arrays parametrizing the Gel$'$fand-Cetlin basis for 
  representations of $U(n,\H)$ (or $Sp(2n,\C)$).
  The top row \(\lambda = (\lambda_1, \ldots,
  \lambda_n)\) is fixed, and is given by the highest weight of the
  irreducible $V(\lambda)$.} \label{fig:GCforUnH_2}
}
\end{center}
\end{figure}

Molev finds another approach in \cite{M4}.
There is an algebra map
\begin{equation}\label{eq:MolevmapPsi}
\Psi: Y^{-}(2) \to U(\sp(2n,\C))^{\sp(2(n-1),\C)},
\end{equation}
where $Y^{-}(2)$ is an infinite-dimensional algebra called the {\em
twisted Yangian}.  Molev then shows that the induced action of
$Y^{-}(2)$ on $(V(\lambda)^{\mu})^{+}$ 
is still irreducible. This map $\Psi$, originally used by
Ol$'$shanskii in \cite{O92} and simplified by Molev and Ol$'$shanskii in
\cite{MO}, is the key new ingredient to Molev's construction of a
Gel$'$fand-Cetlin basis for $Sp(2n,\C)$. This is because the
representations of Yangians and twisted Yangians are
well-understood, and a Gel$'$fand-Cetlin-type basis for
representations of Yangians is constructed in \cite{M1}.  Molev
explicitly identifies $(V(\lambda)^{\mu})^{+}$, and therefore $\mult$,
with known representations of the Yangian. He then combines the Thimm
action on $V(\lambda)$ with the Yangian action on the multiplicity
spaces to construct a Gel$'$fand-Cetlin basis for representations of
$Sp(2n,\C)$. He finds that, as in the $U(n,\C)$ case, the basis
vectors are parametrized by patterns of integer arrays as in Figure
\ref{fig:GCforUnH_2}.  The fact that the decompositions into
  irreducible $V(\mu)$ under the action of $Sp(2(n-1),\C)$ are {\em
    not} necessarily multiplicity-free is reflected by the presence of the
  additional integer parameters $\lambda' = ({\lambda}'_1, \ldots,
  {\lambda}'_n)$ ``in between'' the
  $\lambda=(\lambda_1,\ldots,\lambda_n)$ and $\mu=(\mu_1,\ldots,
  \mu_{n-1})$, et cetera.

\section{The Gel$'$fand-Cetlin-Molev integrable system}\label{sec:GCM}

In this section we construct the Gel$'$fand-Cetlin-Molev integrable
system on coadjoint orbits $\Ol$ of $U(n,\H)$, which will be the
$U(n,\H)$-analogue of the Gel$'$fand-Cetlin system described in
Section~\ref{subsec:GCSystemforUnC}, and is the answer to the Question
posed in the Introduction. In the Guillemin-Sternberg construction of
the Gel$'$fand-Cetlin system on $U(n,\C)$ coadjoint orbits, the Thimm
functions, obtained by projections to smaller and smaller Weyl
chambers, provide enough functionally independent Poisson-commuting
functions to produce a half-dimensional torus action on the coadjoint
orbits. A simple dimension count reveals that, in the case of
$U(n,\H)$-coadjoint orbits, this is simply not possible: the rank of
the maximal torus in $U(n,\H)$ is too small in comparison to the
dimension of the group.  This problem is the symplectic-geometric
manifestation of the fact that the decomposition of $V(\lambda)$ by
$Sp(2(n-1),\C)$ is not multiplicity-free. See the bottom line in the
``Rosetta Stone'' in Section~\ref{subsec:GCSystemforUnC}.

We construct the Gel$'$fand-Cetlin-Molev system on $\Ol$ in
Section~\ref{subsec:main}. We will explain the interpretations in
terms of the non-trivial reduced spaces $\Ored$ in
Section~\ref{subsec:interpret}.

\subsection{The construction}\label{subsec:main}

We refer the reader to Appendix A for reminders on linear
algebra over $\H$ and Lie-group-theoretic facts about $U(n,\H)$.

The first part of the construction of an integrable system on a
coadjoint orbit $\Ol$
of $U(n,\H)$ follows exactly the procedure used
to construct the Gel$'$fand-Cetlin system on orbits of
$U(n,\C)$. Again, we choose a chain of subgroups
\begin{equation}\label{eq:UnHchain}
U(1,\H) \subset U(2,\H) \subset \cdots \subset U(n-1,\H) \subset
U(n,\H),
\end{equation}
where $U(k,\H)$ is the subgroup of upper left \(k \times k\) matrices in
$U(n,\H)$. 
Recall that the Thimm functions for the Gel$'$fand-Cetlin system 
are obtained by taking projections at
each step to the positive Weyl chamber. The same method works for the
$U(n,\H)$ case to produce \(n(n-1)/2\) Poisson-commuting, independent
functions on $\Ol$. We have the diagram
\begin{equation}\label{eq:ThimmforUnH}
\xymatrix{
\Ol  \ar @{^{(}->}[r] & \unH^{*} \ar[r]  & 
\unHs^{*} \ar[r] \ar[d]^{/U(n-1,\H)} & \u(n-2,\H)^{*} \ar[r] \ar[d]^{/U(n-2,\H)} & \cdots \\
   &  & (\t^{n-1})_{+}^{*}   & (\t^{n-2})^{*}_{+}  & \cdots 
}
\end{equation}
analogous to~\eqref{eq:ThimmforUnC}. In linear-algebraic terms, these
Thimm functions are obtained by diagonalizing a matrix in $\u(k,\H)^*
\cong \u(k,\H)$ by an element of $U(n,\H)$ to a diagonal matrix of the
form~\eqref{eq:UnHdiagonal}, and reading off the diagonal entries
(ignoring factors of $i$). Again, the first projection from \(\unH^*
\to (\t^n)^*_{+}\) is omitted since it is trivial when restricted to
$\Ol$. 

The main obstacle in the $U(n,\H)$ case is that these $n(n-1)/2$
functions do
not suffice to completely integrate a generic coadjoint orbit of
$U(n,\H)$, since for such an $\Ol$ we have
\[
\mathrm{dim}(\Ol) = \mathrm{dim}(U(n,\H)) - \mathrm{dim}(T) = 2n^2.
\]
Thus, to completely integrate a generic $\Ol$, it is necessary to find
an {\em additional} \(n^2 - n(n-1)/2 = n(n+1)/2\) independent, Poisson-commuting
functions on $\Ol$. This will be our main task in this section.

It will turn out that these new functions are obtained by
augmenting the diagram~\eqref{eq:ThimmforUnH} with new functions $G_{n-k}$
as follows.
\begin{equation}\label{eq:ThimmplusGs}
\xymatrix{ & \R^n & \R^{n-1} & \R^{n-2} & \cdots \\ 
     \Ol \ar
     @{^{(}->}[r] & \unH^{*} \ar[r]  \ar @{.>}[u]_{G_n} &
     \unHs^{*} \ar[r] \ar[d]^{/U(n-1,\H)} \ar @{.>}[u]_{G_{n-1}} &
     \u(n-2,\H)^{*} \ar[r] \ar[d]^{/U(n-2,\H)} \ar @{.>}[u]_{G_{n-2}} &
     \cdots \\ 
     &  & (\t^{n-1})_{+}^{*} &
     (\t^{n-2})^{*}_{+} & \cdots }
\end{equation}
Since each \(G_{n-k}\) has $n-k$ components, this is exactly the
number that we need to completely integrate $\Ol$. 
The motivation behind the definitions of these functions $G_{n-k}$ is
the subject of Section~\ref{sec:limits}. In this section, we will
simply take the $G_{n-k}$ as defined, and 
concentrate on showing that they (along with the Thimm functions)
integrate $\Ol$.

We will now define these new functions $G_{n-k}$. Since they are all
  defined analogously, for concreteness we define $G_n$. 
  We denote the components of $G_n$ by \(G_n = (g_{n,1}, g_{n,2}, \ldots,
  g_{n, n}).\) We begin with the first $n-1$ components \((g_{n,1}, \ldots,
g_{n,n-1}).\) Let \(X \in
\unH^* \cong \unH\) be an element in the $U(n,\H)$-orbit of
$\intwc$. Then there exists a unique diagonal matrix \(D_{\lambda} :=
(i\lambda_1,\ldots,i\lambda_n),\) where 
\[
0 > \lambda_1 > \lambda_2 > \ldots > \lambda_n,
\]
such that $X$ is conjugate by $U(n,\H)$ to $D_{\lambda}$. Let 
\[X = A  D_{\lambda}  A^*\] 
for an element \(A \in U(n,\H).\) Note this equation defines $A$ up to
right multiplication by the maximal torus $T^n$ of $U(n,\H)$. 
We take 
\begin{equation}\label{eq:gmdef}
g_{n,m}(X) := |a_{n,m}|^2, \quad 1 \leq m \leq n-1. 
\end{equation}
So the $g_{n,m}$ just takes a norm-square of an entry in the bottom
row of the matrix $A$.  Since the norm-squares are $T^n$-invariant,
the $g_{n,m}$ are well-defined. 

\begin{remark}\label{remark:Gnn}
From this description, it is clear that we cannot define the $n$-th
component $g_{n,n}$ of $G_n$ in the same way as the first $n-1$
components, because (since $A$ is unitary) the matrix entries always satisfy
\[
\sum_{m=1}^n |a_{n,m}|^2 = 1.
\]
Thus the components in $G_n$ would {\em not} be functionally
independent if $g_{n,n}(X)$ were also defined to be $|a_{n,n}|^2$. 
\end{remark}

We now define the $n$-th component of $G_n$, which has a qualitatively
different description. Recall that the maximal
torus $T^n$ also acts on the coadjoint orbit $\Ol$, with moment map
induced by the inclusion \(\iota: \t^n \into \unH:\) 
\begin{equation}\label{eq:Tmommap}
\xymatrix{
\Ol \ar @{^{(}->}[r] & 
\unH^* \ar[r]^{\iota^*} & (\t^n)^* \cong \R^n.
}
\end{equation}
Here we take the standard identification of $\t^n$ with its dual to
identify $(\t^n)^*$ with $\R^n$. We define the $n$-th component
$g_{n,n}$ of $G_n$ to be the $n$-th component of the moment
map~\eqref{eq:Tmommap}, i.e. for \(\epsilon_n\) the standard $n$-th
basis vector in \(\R^n \cong \t^n,\) we have 
\begin{equation}\label{eq:gnn}
g_{n,n}(X) := \langle \iota^*(X), \epsilon_n \rangle. 
\end{equation} 
Note that in the $U(n,\H)$ case, the components of the moment map for
the action of the maximal torus are functionally independent of the
Thimm trick functions, in contrast to the $U(n,\C)$ case. Hence it
makes sense to use them as components of the $G_{n-k}$. 

The functions $G_{n-k}$, as mentioned above, are defined
analogously. Before stating the main results, we make a remark on
notation. 
From the sequence of 
subgroups~\eqref{eq:UnHchain}, we get a sequence of moment maps 
\begin{equation}\label{eq:PhikUnH}
\xymatrix{
\Ol \ar @{^{(}->}[r]^(0.4){\Phi} \ar @/^2pc/[rr]_{\Phi_{n-1}} \ar @/^3pc/[rrr]^{\Phi_{n-2}} 
& \unH^{*} \ar[r] & 
 \unHs^{*} \ar[r] &  \u(n-2,\H)^{*} \ar[r] &  \cdots  
}
\end{equation}
By abuse of notation, we will sometimes denote the pullbacks
\(\Phi_{n-k}^*G_{n-k}\) by $G_{n-k}$. Similarly, we will sometimes
refer to the Thimm functions \(\u(n-k,\H)^* \to (\t^{n-k})^*_{+}\) as
functions on $\Ol$, by pulling back via $\Phi_{n-k}$.  With this
notation in place, the main theorem of this paper may now be stated.

\begin{theorem}\label{thm:main}
Let $\Ol$ be a generic coadjoint orbit of $U(n,\H)$.  The functions
$\{G_{n-k}\}_{k=0}^{n-1}$ as defined in~\eqref{eq:gmdef}
and~\eqref{eq:gnn}, plus the Thimm functions defined
in~\eqref{eq:ThimmforUnH}, pull back via the
diagram~\eqref{eq:ThimmplusGs} to give a completely integrable system
on an open dense set of $\Ol$.
\end{theorem}

We call this the {\em Gel$'$fand-Cetlin-Molev system} on the coadjoint orbit
$\Ol$ of $U(n,\H)$. To prove Theorem~\ref{thm:main}, 
there are two things to check: that the functions above
Poisson-commute, and that they are functionally independent. Thus
Theorem~\ref{thm:main} follows immediately from the following two
propositions. 

\begin{proposition}\label{prop:Pcommute}
Let $\Ol$ be a generic coadjoint orbit of $U(n,\H)$.  Let
$\{G_{n-k}\}_{k=0}^{n-1}$ be defined as in~\eqref{eq:gmdef}
and~\eqref{eq:gnn}, and the Thimm functions defined
in~\eqref{eq:ThimmforUnH}.  Then these functions Poisson-commute on an
open dense subset of $\Ol$.
\end{proposition}

\begin{proposition}\label{prop:indep}
Let $\Ol$ be a generic coadjoint orbit of $U(n,\H)$.  Let
$\{G_{n-k}\}_{k=0}^{n-1}$ be defined as in~\eqref{eq:gmdef}
and~\eqref{eq:gnn}, and the Thimm functions defined as
in~\eqref{eq:ThimmforUnH}.  Then these functions are independent on an
open dense subset of $\Ol$.
\end{proposition}

In the course of the proofs of both Propositions, it will turn out to
be convenient to replace the functions $G_{n-k}$ with functions
$F_{n-k}=(f_{1,n-k}, \ldots, f_{n-k, n-k})$ 
defined as follows. Again, for concreteness we define $F_n$,
but the others are defined similarly. 
\begin{equation}\label{eq:fmn}
f_{n,m}(X) := \sum_{\ell =1}^{n} (-1)^m \lambda_{\ell}^{2m}
|a_{n,\ell}|^2.
\end{equation}
The motivation behind the definitions of $F_{n-k}$ will be explained
fully in Section~\ref{sec:limits}. Indeed, it is the $F_{n-k}$, and not the $G_{n-k}$, which are
obtained directly, in
Theorem~\ref{thm:fmnDerive}, as classical limits of generators of an {\em
abelian} subalgebra in a (non-commutative) algebra $Y^{-}(2)$.

Going back to~\eqref{eq:fmn}, we first note that the functions $f_{n,m}$ are obtained from the
$g_{n,m}$ and from the components of the Thimm function to
$(\t^n)^*_{+}$. Since the matrix 
\[
\left[
\begin{array}{cccc}
(-1) \lambda_1^2 & (-1)^2 \lambda_1^4 & \cdots & (-1)^{n}(\lambda_1)^{2n} \\
(-1) \lambda_2^2 & (-1)^2 \lambda_2^4 &        & (-1)^{n}(\lambda_2)^{2n} \\
\vdots &  & & \vdots \\
(-1)\lambda_n^2 & (-1)^2 \lambda_n^4 & \cdots & (-1)^{n}(\lambda_n)^{2n} \\
\end{array}
\right].
\]
is invertible when the $\lambda_{\ell}$ are distinct, the $g_{n,m}$
may also be obtained from the $f_{n,m}$ plus the Thimm
functions. Thus, for the purposes of showing independence and
Poisson-commutativity, it is equivalent to use the
$F_{n-k}$. 

\begin{remark}
Remark~\ref{remark:Gnn} also applies to the $F_{n-k}$, in that the $n$
functions $\{f_{1,n}, \ldots, f_{n,n}\}$ are not independent, but only
give $n-1$ independent functions. In either case, it is necessary
to also include the ``extra'' $S^1$ moment map as given in
equation~\eqref{eq:gnn}.
\end{remark}

It will also be useful to have in hand another
description of the functions $F_{n-k}$. Again, for simplicity, we take
the case \(k=0.\) Given an element \(X = AD_{\lambda}A^* \in \unH^*
\cong \unH,\) it is a straightforward computation to verify that
\begin{equation}\label{eq:fmn2}
f_{n,m}(X) = \mbox{rtr}(X^{2m}E_{nn}) =
\mbox{rtr}(AD_{\lambda}^{2m}A^*E_{nn})
\end{equation}
is equivalent to the formula given in~\eqref{eq:fmn}. Here
$\mathrm{rtr}$ denotes the reduced trace pairing on $\unH$ defined
in~\eqref{reducedtrace}. 
We will use this form in the proofs below.

\medskip

\begin{proof} [of Proposition~\ref{prop:Pcommute}] 

The symplectic leaves on the dual of a Lie algebra $\g^*$ are 
the orbits under the
coadjoint action of $G$, so any $G$-invariant function on
$\g^*$ is a Casimir.
Since the Thimm function on $\u(n-k,\H)^*$ is by construction
$U(n-k,\H)$-invariant, this implies that the components of the Thimm function from
$\u(n-k,\H)^*$ to \((\t^k)^*_{+}\) Poisson-commutes with any component
of $G_{n-k}$. Similarly, they Poisson-commute with anything ``to the
right'' in the diagram~\eqref{eq:ThimmplusGs}, i.e. any component of
$G_{n-p}$ or of the Thimm functions from $\u(n-p,\H)^*$ for any \(k <
p \leq n-1.\) 
By a similar argument, since the components of $G_{n-k}$ are
$U(n-k-1,\H)$-invariant, they Poisson-commute with any component of
the projection to $\u(n-k-1,\H)^*$. 

It remains to show that, at each step, the components of
$G_{n-k}$ Poisson-commute with each other. For concreteness, we consider the case
$k=0$. The other steps may be argued similarly. 
In fact, as
remarked above, it will here be more convenient to use the function
$F_n$ as defined in~\eqref{eq:fmn2} rather than the $G_{n-k}$. 

We will first show that $g_{n,n}$,
the $S^1$ moment map, commutes with all \(f_{n,m}, 1 \leq m \leq n.\)
Note that it suffices to check that they commute on a fixed symplectic
leaf ${\cal O}$. For the case \(k=0\) which we consider, the only relevant
symplectic leaf is the original coadjoint orbit $\Ol$. 
Here and below, we will use for convenience the projection map \(\pi:
U(n,\H) \to U(n,\H)/T^n \cong {\cal O}\) to pull back calculations to
$U(n,\H)$. Let \(\overline{f}_{m,n}\) denote the pullback \(\pi^*f_{n,m}\) for \(1 \leq m \leq
n.\) 
Let $Y^{\sharp}$ denote the vector field on ${\cal O}$ generated by the
$S^1$-action, and let \(\overline{Y^{\sharp}}\) denote the
corresponding vector field on $U(n,\H)$ generated by the $S^1$ (by
left multiplication). By construction,
\(d\pi(\overline{Y^{\sharp}}) = Y^{\sharp}.\) 
By definition of the Poisson bracket, it suffices to show that
\[
df_{n,m}(Y^{\sharp}) \equiv \{f_{n,m}, g_{n,n}\} \equiv 0,
\]
for any $1 \leq m\leq n$. 
Lifting to $U(n,\H)$, it suffices to show that 
\[
d\overline{f}_{m,n}(\overline{Y^{\sharp}}) = 0. 
\]
Let \(A \in U(n,\H).\) We trivialize the tangent bundle to $U(n,\H)$
by right multiplication. Let \(W^{\sharp}_A\) denote the tangent vector
at $A$ corresponding to \(W \in \unH \cong T_{1}(U(n,\H)).\) Using the
definition of $\overline{f}_{m,n}$, it is straightforward to compute that 
\begin{equation}\label{eq:dgm}
d\overline{f}_{m,n}(W^{\sharp}_A) = - \text{rtr}([(A D_{\lambda} A^*)^{2k},E_{nn}]W).
\end{equation}
Since the Hamiltonian vector field $Y^{\sharp}$ associated to the $S^1$-action
is generated by the element \(Y = iE_{nn} \in \unH,\) 
we immediately find that 
\[
d\overline{f}_{m,n}(\overline{Y^{\sharp}}) = 0. 
\]
Therefore, by definition, $g_{n,n}$ Poisson-commutes with any $f_{n,m},1 \leq m \leq n$. 

Finally, it remains to show that the $f_{n,m}$, for \(1 \leq m \leq n,\)
Poisson-commute among themselves. This fact follows from the construction
of the $f_{n,m}$ as classical limits of generators of an abelian
subalgebra contained in $Y^{-}(2)$. Since the corresponding generators
commute in the quantization, the classical limits automatically {\em
  Poisson-}commute. 

\end{proof}

\begin{remark}
It is possible to prove directly, using the standard Kostant-Kirillov
Poisson structure on $\unH^*$, that the $f_{n,m}$
Poisson-commute. However, the calculation is long, and we thus
prefer to invoke the deformation theory. 
\end{remark}

We must now show that the functions in~\eqref{eq:ThimmplusGs} are independent. 
In fact, we will show more: 
they induce independent functions on the
reduced spaces.This interpretation will be further discussed in
Section~\ref{subsec:interpret}.

\medskip
\begin{proof}[of Proposition~\ref{prop:indep}]

As in the proof of Proposition~\ref{prop:Pcommute}, we will
occasionally use the $F_{n-k}$ instead of the $G_{n-k}$. 
We first claim that the $\{f_{n,m}\}_{m=1}^{n-1}$ are
independent. This is clear, since they are norm-squares of different
matrix entries of elements in $U(n,\H)$. Similarly, the components
\(\{f_{n-k,m}\}_{m=1}^{n-k-1}\) are also independent for any \(1 \leq k \leq n-1.\)  
Second, we claim the
components of $G_{n-k}$ for a fixed $k$, \(1 \leq k \leq n-1\) are
also independent of the components of the Thimm function
\(\u(n-k,\H)^* \to (\t^{n-k})^*_{+}.\) This is because on the open
dense set \(U(n-k,\H) \cdot ((\t^{n-k})^*_{+})_{0} \cong
U(n-k,\H)/T^{n-k} \times ((\t^{n-k})^*_{+})_{0} \subset
\u(n-k,\H)^*,\) the Thimm function simply reads off the second factor,
whereas the components of $G_{n-k}$ are functions on the first factor.
Third, for a fixed $k$, \(1 \leq k \leq n-1,\) the
last component $g_{n-k,n-k}$ is also independent of the Thimm
functions. This is because $g_{n-k,n-k}$ generates, by
definition, a non-trivial action on each generic symplectic leaf in
$\u(n-k,\H)^*$, and in particular is non-constant on those leaves,
whereas the
Thimm functions are constant on leaves. 
Fourth, the components of the Thimm function
\(\u(n-k,\H)^* \to (\t^{n-k})^*_{+}\) are independent of any component
of the projection \(\iota^*: \u(n-k,\H)^* \to \u(n-k-1,\H)^*.\) To see this,
note that the restriction of $\iota^*$ to a coadjoint orbit ${\cal O}$
in $\u(n-k,\H)^*$ gives the moment map $\mu_{\cal O}$ for the $U(n-k-1,\H)$-action on
${\cal O}$. We assume ${\cal O}$ is an orbit through a point $\nu$ in the
interior of $(\t^{n-k})^*_{+}$. 
For any \(X \in \u(n-k-1,\H),\)
there exists an element \(Z \in \u(n-k,\H)\) such that 
\[
d\mu_{\cal O}^X(Z^{\sharp}) = \omega_{\nu}(X^{\sharp}, Z^{\sharp}) =
\left<\nu, [X,Z]\right> \neq 0.
\]
This implies that $\mu_{\cal O}^X$ is non-constant on ${\cal O}$ for
all \(X \in \u(n-k,\H).\) In particular, any component is functionally
independent of the Thimm functions (which are constant on ${\cal
  O}$).

It remains now to show (and this is the bulk of the proof) two things: 
that the last component $g_{n-k,n-k}$ of $G_{n-k}$ 
is independent of the first $n-k-1$ components,
and that all components of $G_{n-k}$ are independent of any component
of the projection \(\u(n-k,\H)^* \to \u(n-k-1,\H)^*.\) Without loss of
generality we consider the case \(k=0.\) 

In order to show that the last component $g_{n,n}$ is independent of the
first $n-1$ components of $G_n$, 
it suffices to show that there exists
an element \(Z \in \unH^*\) such that on $\Ol$ we have 
\[
d{f_{n,m}}(Z^{\sharp}) \equiv 0
\]
for all \(1 \leq m \leq n-1,\) but 
\[
dg_{n,n}(Z^{\sharp}) \not \equiv 0.
\]
Consider the subgroup \(U(1,\H)\) sitting in $U(n,\H)$ as the ``bottom
right'' $(n,n)$-th entry. For any \(v \in
\text{Im}(\H),\) where \(a,b,c \in \R\) and \(\|v\| = 1,\) there is a
corresponding $S^1$ subgroup
\(\{e^{i\theta} \vert \theta \in \R\} \subset U(1,\H).\) 
Let 
$Z_v^{\sharp}$ denote the vector field generated on $\Ol$ by this copy of
$S^1$. Then at a point \(\xi \in \Ol,\) 
\begin{eqnarray*}
d(g_{n,n})_{\xi}(Z_v^{\sharp}) &=& \omega_{\xi}(Y^{\sharp}, Z_v^{\sharp}) \\
 &  = & \langle \xi, [Y,Z] \rangle.
\end{eqnarray*}
Here $Y =iE_{n,n}$ is the element in $\unH$ generating the
$S^1$-action corresponding to $g_{n,n}$, and $Y^{\sharp}$ is the
corresponding vector field on $\Ol$. 
By using the action of $U(n-1,\H)$, we may 
take the point $\xi \in \unH^*$ to be of the form
\begin{equation}\label{eq:normalize}
\xi = \left[\begin{array}{cc} D_{\mu} & \star \\
                            \star & z \end{array} \right],
\end{equation}
for \(z \in \text{Im}(\H),\) and a diagonal matrix $D_{\mu}$ 
of the form~\eqref{eq:UnHdiagonal}. 
Since \(\xi \in \Ol,\) the $(n,n)$-th entry $z$ cannot be equal to
zero. In particular, there exists $v \in \text{Im}(\H),\|v\| =1,$ for which
\(Z = v E_{n,n} \in \unH\) has the property that 
\begin{eqnarray*}
d(g_{n,n})_{\xi}(Z_v^{\sharp}) & = & \langle  \xi, [Y,Z] \rangle \\
 & = & \text{rtr}(\xi[Y,Z]) \\
 & = & \mathrm{Re}(z(i \cdot v - v \cdot i)) \in \H \\
 & \not = & 0.
\end{eqnarray*}
On the other hand, from equation~\eqref{eq:dgm}, we conclude that 
for \(v \in \text{Im}(\H),\) we have that 
\[
d{f_{n,m}}(Z_v^{\sharp}) = 0
\]
for \(1 \leq m \leq n-1.\) Thus $g_{n,n}$ is independent from the
$\{f_{n,m}\}_{m=1}^{n-1}$. 
In fact, this argument additionally shows that $g_{n,n}$ is
independent of any component of the moment map \(\Phi_{n-1}: \Ol \to
\u(n-1,\H)^*,\) since by construction, \([W, Z_v]=0\) for any \(W \in
\u(n-1,\H).\)

Our last task is to show that the components of $F_n$ are
independent from the components of the moment map \(\Phi_{n-1}: \Ol \to
\u(n-1,\H)^*.\) As in
the proof of Proposition~\ref{prop:Pcommute}, we will do 
calculations on $U(n,\H)$ instead of $U(n,\H)/T$. Define
\(\overline{\Phi}_{n-1} := \pi^*(\Phi_{n-1})\) and \(\overline{F}_n :=
\pi^*(F_n).\) 
It will suffice to show that there exists some $A \in U(n,\H)$
such that 
the linear equations defining the kernels of both $d\overline{\Phi}_{n-1}$
and $d\overline{F}_n$ are linearly independent at $A$.

Let \(A \in U(n,\H).\) We trivialize $TU(n,\H)$ by right
translation. Let \(X^{\sharp}_{A}\) denote the right translate of
\(X\)
to $T_AU(n,\H)$. Let $\jmath$ denote the map \(\u(n,\H) \to
\u(n-1,\H)\) given by taking the upper left \((n-1)\times (n-1)\)
submatrix.  Then the pullback of the moment map
\(\overline{\Phi}_{n-1}\) can be expressed as
\[
\overline{\Phi}_{n-1}: A  \mapsto \jmath(A D_{\lambda} A^{*}),
\]
and the derivative by 
\begin{equation}\label{eq:kerPullbackPhi}
(d\overline{\Phi}_{n-1})_{A}(X^{\sharp}_A) 
  =   \jmath([X, AD_{\lambda} A^{*}]).
\end{equation}
Since $\jmath$ is the map which takes the upper left submatrix,
it is convenient to write an element \(X \in \unH\) as
\begin{equation}\label{eq:Xdecomp}
X = \left[ \begin{array}{cc} X_{11} & X_{12} \\ - X_{12}^{*} & X_{22} 
\end{array} \right],
\end{equation}
where \(X_{11} \in \unHs, X_{12} \in \H^{n-1}, X_{22} \in \text{Im}(\H).\) 
Moreover, as in~\eqref{eq:normalize}, we may 
assume that $A$ is such that 
\[
A D_{\lambda} A^{*} = \left[\begin{array}{cc} D_{\mu} & W \\ -W^{*} & z
\end{array} \right],
\]
for \(W \in \H^{n-1}\) and $z$ nonzero. 
Then from~\eqref{eq:kerPullbackPhi} we see that the linear equations (over $\R$) 
defining \(\kerPhi \in T_{A}(U(n,\H)) \cong \unH\) are given by
the single matrix equation
\[
X_{11} D_{\mu} - X_{12} W^{*} = D_{\mu} X_{11} + W X_{21}.
\]
In other words, the matrix \(X_{11} D_{\mu} - X_{12} W^{*}\)
must be $\H$-hermitian. 

Now we compute the linear equations for \(\text{ker}(\overline{f}_{m,n}).\) 
The pullbacks $\overline{f}_{m,n}$ are given by 
\begin{equation}\label{eq:barfmdef}
\overline{f}_{m,n}: \quad A \mapsto -\mathrm{rtr}(A D_{\lambda}^{2m}
A^{*} E_{nn}). 
\end{equation}
Let $X$ be written as in~\eqref{eq:Xdecomp}, and write
\(X_{12} = ((X_{12})_1, \ldots, (X_{12})_{n-1})^t \in \H^{n-1},\)
where \((X_{12})_i \in \H.\) Then 
\begin{eqnarray*}\label{eq:dbarfm}
d\overline{f}_{m,n}(X^{\sharp}_{A}) & = & -\mathrm{rtr} \left([A
 D_{\lambda}^{2m} A^{*}, E_{nn}] X \right) \\ 
 & = & -2 \cdot
 \mathrm{Re} \left(\sum_{i=1}^{n-1} (-1)^m \left(\sum_{\ell=1}^n
 \lambda_\ell^{2m} a_{i\ell}\overline{a}_{n\ell} \right)(- \overline{(X_{12})_i}) \right)
 \\ 
 & & + 2 \cdot \mathrm{Re} \left(\sum_{i=1}^{n-1}(-1)^m
 \left(\sum_{\ell=1}^n \lambda_\ell^{2k} a_{n\ell}\overline{a}_{i\ell} \right)
 (X_{12})_{i} \right),
\end{eqnarray*}
so \(\mathrm{ker}d\overline{f}_{m,n}\) is given by the condition that
the expression above is $0$.  Denote the (real)
variables in $X_{11}$ by $\{z_a\}$, and the (real) variables in
$X_{12}$ by $\{w_b\}$.  Note that the linear equations defining the
$\text{ker}(\overline{f}_{m,n})$, as seen above, involve only the
$w_b$.  The linear equations defining $\kerPhi$ involve both $z_a$ and
$w_b$, but they are also linearly independent modulo $\langle w_b
\rangle$.  
This implies that the \(d\overline{f}_{m,n}\) are
also linearly independent on \(\kerPhi.\) Thus the \(\{f_{n,m}\}\) are
independent of any component of $\Phi_{n-1}$, as desired. 

\end{proof}

\subsection{Intepretation on the reduced
  spaces}\label{subsec:interpret}

We now take a moment to interpret the results of the previous section in terms of the
reduced spaces $\Ored$, and comment on the differences between the
cases of $U(n,\C)$ and $U(n,\H)$. 

We already mentioned in the beginning of Section~\ref{subsec:main}
that the essential new problem in the $U(n,\H)$ case is that the Thimm
functions do not give ``enough'' functions to completely integrate a
generic coadjoint orbit of $U(n,\H)$. We now briefly review the Thimm
trick construction of torus actions, and explicitly see that this deficiency
is due to the presence of non-trivial symplectic reductions. 

Suppose a compact Lie group $G$ acts on a symplectic manifold $M$ with
moment map $\mu$. Let \(T \subset G\) be a maximal torus, and identify
\(\g \cong \g^*\) so that \(\t^* \subset \g^*.\) We assume $\mu(M)
\cap \intwc \neq \emptyset$. It is shown in
\cite{LMTW} that the preimage \(S := \mu^{-1}(\intwc)\) is a
symplectic submanifold of $M$, and the restriction of $\mu$ to $S$ is
a moment map for the $T$-action on $S$.  This submanifold $S$ is
called a {\em symplectic slice.} Note that for a regular value \(\alpha
\in \intwc,\) we have \(M \mmod_{\alpha} G = S \mmod_{\alpha} T.\) The
Thimm torus is then defined to act on $G \cdot S$ as follows: for any
\(g \cdot p \in G \cdot S\) and \(t \in T,\) define 
\[
t \cdot (g \cdot p) := g \cdot (t \cdot p),
\]
where by \(t \cdot p\) we mean the {\em original} $T$-action on $M$,
restricted to the slice $S$. 

Suppose that the symplectic manifold is a coadjoint orbit of $U(n,\C)$
or $U(n,\H)$ and $G$ is the first subgroup in the
chain~\eqref{eq:chain} or~\eqref{eq:UnHchain}, respectively. Let
\(p \in S.\) (We consider points in the slice without loss of
generality; for any other point in \(G \cdot S\) we could repeat the
argument with a conjugate torus.) The presence of a
completely integrable system translates to the presence of a
half-dimensional Lagrangian subspace \(L_p \subset T_pM\) spanned by
the Hamiltonian vector fields of the Poisson-commuting functions. By
the above description of the Thimm torus action, we see that the
Hamiltonian vector fields arising from the Thimm functions are exactly
the $X^{\sharp}_p$ for \(X \in \t \subset \g,\) so the span is exactly
\(T_p(T \cdot p) \subset T_pS,\) giving an isotropic subspace of
$T_pS$. There are two reasons why $T_p(T \cdot p)$ may not be a
Lagrangian subspace of $T_pM$. 
First, perhaps most obviously, $T_pS$
is not all of $T_pM$.  There is a complementary symplectic subspace of
$T_pS$ in $T_pM$, mapping isomorphically under $d\mu$ to the tangent
space to the coadjoint orbit ${\cal O}_{\mu(p)}$, which is not
accounted for by $T_p(T \cdot p)$. Second, $T_p(T \cdot p)$ may not be
Lagrangian even in $T_pS$. There is a symplectic subspace in $T_pS$
mapping isomorphically to the tangent space of the reduced space
$T_{[p]}(M \mmod T)$ which is also not accounted for by the $T_p(T
\cdot p)$. 

The first reason mentioned above, the presence of a subspace
isomorphic to $T_{\mu(p)}{\cal O}_{\mu(p)}$, is partially resolved by
the inductive step. Namely, ${\cal O}_{\mu(p)}$ is itself a symplectic
manifold with respect to $H$, where $H$ is now either \(U(n-2,\C\) or
\(U(n-2,\H).\) 
By considering the
Thimm functions arising from the action of $H$, the subspace for
$T{\cal O}_{\mu(p)}$ will in turn break up into pieces, part of which
will be spanned by the Hamiltonian vector fields arising from the
Thimm functions from $H$. 
Hence it is the second reason mentioned above which is the
essential obstacle to having the Thimm functions completely integrate
the original manifold $M = \Ol$. At each inductive step, if the
subspace corresponding to the reduced space is
non-trivial and thus has positive dimension, then it is impossible for
the Thimm functions to integrate $\Ol$. 

We may now compare the cases of $U(n,\C)$ and $U(n,\H)$ in this
light. The Gel$'$fand-Cetlin construction given in
Section~\ref{subsec:GCSystemforUnC} works exactly because the
symplectic reductions ${\cal O} \mmod U(n-1,\C)$ are {\em trivial}. The
analogous construction for the $U(n,\H)$ case does not work, and we
need more functions, precisely because the symplectic reductions
\({\cal O} \mmod U(n-1,\H)\) are {\em not} trivial. Indeed, generically
they are dimension $2n$. This is in exact correspondence with the
positive-dimensionality of the multiplicity spaces $M_{\lambda}^{\mu}$
in Section~\ref{subsec:GCforUnH}. We invite the reader to take another
look at the Rosetta Stone in Section~\ref{subsec:GCforUnC} with these
interpretations in mind. 

As advertised in the previous section, the functions $G_n$ may be viewed as an integrable system on the reduced spaces. We record the following, which we already stated in the Introduction as Theorem~\ref{thm:Oredfunctions}.

\begin{theorem}\label{thm:Oredsystem}
Let \({\mathcal O}_{\lambda} \cong U(n,\H)/T^n\) be a coadjoint orbit of 
$U(n,\H)$. Let $G_n$ be defined as in equations~\eqref{eq:gmdef}
and~\eqref{eq:gnn}. 
Then the components of $G_n$ 
descend to functionally independent, Poisson-commuting functions on
the reduced space \({\mathcal O}_{\lambda} \mmod\!_{\mu} U(n-1,\H)\)
for $\mu$ a regular value.  
\end{theorem} 

\begin{proof}
Since the $G_n$ defined by~\eqref{eq:gmdef}
and~\eqref{eq:gnn} are $U(n-1,\H)$-invariant, and because they are
shown in the proof of Proposition~\ref{prop:indep} to be functionally
indepedent of any component of $\Phi_{n-1}$, they automatically
induce functionally independent functions on the reduced spaces \({\cal
  O} \mmod U(n-1,\H).\) Moreover, by the definition of the
symplectic structure on the reduced space, they also automatically
Poisson-commute on the reduced space. Thus we have $n$ independent
Poisson-commuting functions, and therefore a completely integrable
system, on the reduced space.
\end{proof}

\section{The classical limits}\label{sec:limits}

In this section, we explain our deformation-theoretic derivation of the formul{\ae} for the
functions $f_{n,m}$ used in the construction of the Gel$'$fand-Cetlin-Molev integrable system.
Some comments are in order: first of
all, as seen in the previous section, once the formul{\ae} are given, it is 
possible to prove directly that the functions $f_{n,m}$ give a
completely integrable system on $\Ored$. However, the way in
which these formul{\ae} were derived, via the use of the theory of
deformation quantization and classical limits, is a beautiful story in
its own right. Hence we present it in this section using this
perspective. Since
the construction is the same at each step in~\eqref{eq:ThimmplusGs},
we concentrate here solely on the first step, i.e. the derivation of
the $\{f_{n,m}\}_{m=1}^n$. 

We now recall briefly the basic general philosophy underlying the
computations below. We refer the reader to \cite{CP} for details. The central theme is that the classical limit of 
a non-commutative algebra $A_h$ is a commutative algebra $A_0$
equipped with a Poisson bracket. The Poisson bracket is the
``first-order term'' in the parameter $h$, and the commutative algebra
$A_0$ can then viewed as a space of functions $Fun(M)$ on a Poisson
space $M$. Following
this general recipe, our task in this section is as follows.  We will
first determine in Section~\ref{subsec:CLalgebras} the classical
limits of the non-commutative algebras involved in the Molev
construction as recounted in Section~\ref{subsec:GCforUnH}. In
Section~\ref{subsec:CLmap} we take the classical limit of the algebra
map $\Psi$ used in Section~\ref{subsec:GCforUnH} and get a map between
Poisson spaces. The technical heart of the calculation lies in
Theorem~\ref{thm:basicLimits}, which allows us to obtain  in Theorem~\ref{thm:fmnDerive} the explicit
formul{\ae} for the $f_{n,m}$. 

Some basic necessary definitions and constructions regarding Yangians,
twisted Yangians, and their deformations are recounted briefly in
Appendix B.

\subsection{The classical limits of the algebras}\label{subsec:CLalgebras}

In order to describe the classical limit of the Yangian, following the standard methods in deformation quantization,
we must first exhibit the Yangian as a deformation.  There is a family
of topological Hopf algebras $Y_h(2n)$ for \(h \in \C \backslash
\{0\},\) where \(Y_{1}(2n) = Y(2n),\) and in fact \(Y_h(2n) \cong
Y(2n)\) for \(h \not = 0.\) We then set \(h=0\) in the formul{\ae} for
the algebra and coalgebra structures for $Y_h(2n)$ to obtain a
classical limit.

We define $Y_h(2n)$ as follows. We denote the generators of 
a fixed \(h \in \C \backslash \{0\}\) by \(\tilde{t}_{ij}^{(M)},\)
for \(i, j\) in ${\cal I}$ as in~\eqref{eq:indexingSet} and \(M \geq 0.\) We define the algebra structure
by 
\begin{equation}\label{eq:defYangcommrelations}
[\tilde{t}_{ij}^{(M)}, \tilde{t}_{kl}^{(N)}] = h  \cdot \left(
\sum_{r=0}^{\min(M,L)-1} 
(\tilde{t}_{kj}^{(r)} \tilde{t}_{il}^{(M+L-1-r)} - 
\tilde{t}_{kj}^{(M+L-1-r)} \tilde{t}_{il}^{(r)})\right).
\end{equation}
Note that this is the same formula as for $Y(2n)$ except that we
multiply by a factor of $h$.  The coproduct, antipode, and co-identity
structures are defined by the same formul{\ae} as for $Y(2n)$ with {\em
no} additional factor of $h$.  One can check that these definitions
give $Y_{h}(2n)$ the structure of a Hopf algebra.  We have, by
definition, \(Y_1(2n) = Y(2n).\) In fact, for \(h \not = 0,\)
$Y_h(2n)$ is isomorphic to $Y(2n)$ as a Hopf algebra.

\begin{lemma}\label{lemma:deg_y}
Let $Y_h(2n)$ be defined as above. Then for \(h \not = 0,\) 
\[
Y_h(2n) \cong Y(2n),
\]
as Hopf algebras. 
\end{lemma}

\begin{proof}
The map \(\gamma_h: Y_h(2n) \to Y(2n)\) is given by \(\tilde{t}_{ij}^{(M)}
\mapsto t_{ij}^{(M)}h^M,\) extended linearly. 
\end{proof}

In order to describe the classical limit of $Y(2n)$, we need first
some terminology. Let \(U\) denote a formal
neighborhood of $\infty$ in \(\P^1 = \C \cup \{\infty\}.\) If $u$ is
the usual coordinate on $\C$, then $u^{-1}$ is a coordinate on a
neighborhood of $\infty$.  By a ``formal'' neighborhood, 
we mean that the space of
functions on $U$ is the space of formal power series in the
local coordinate $u^{-1}$.
Let \(U \times \C^{2n}\) be a trivial vector bundle over the formal
neighborhood of $\infty$. A gauge transformation of this vector bundle
is given by an element \(F(u) \in Maps(U, GL(2n,\C)).\) Here, $F(u)$ is a {\em formal}
power series with coefficients in \(\gl(2n,\C)\), with the additional
restriction that it is invertible as a formal power series. Multiplication in the
gauge group is given pointwise. 
The {\em pointed} gauge group \(Maps_{1}(U,GL(2n,\C))\)
is the subgroup such that the point $\infty \in U$ maps
to the identity element in $GL(2n,\C)$. 
Thus an element \(F(u) \in Maps_{1}(U,GL(2n,\C))\) is of the
form 
\[
F(u) = \Id + A_1 u^{-1} + A_2 u^{-2} + \ldots,
\]
where $A_i \in \gl(2n,\C)$. Note that since the $0$-th coefficient is the identity
matrix, such formal power series are invertible for any choice
of \(A_i \in \gl(2n,\C).\) 

\begin{theorem}\label{pl_y}
The classical limit of the Yangian $Y(2n)$ is the 
(infinite-dimensional) pointed gauge group ${\cal G}_{2n} := Maps_{1}(U, GL(2n,\C))$ of
a trivial $\C^{2n}$-bundle over a formal neighborhood $U$ of \(\infty \in
\P^1.\) Moreover, ${\cal G}_{2n}$ has a Poisson structure compatible
with the product structure in the gauge group, making it a Poisson-Lie
group. 
\end{theorem}

\begin{proof}
We denote by $z_{ij}^{(M)}$ the coordinate function on $Maps_{1}(U,
GL(2n,\C))$ which reads off the $(i,j)$-th entry of the coefficient of
$u^{-M}$ of an element \(A(u) = \sum_{M=0}^{\infty} A_M u^{-M} \in
Maps_{1}(U,GL(2n,\C)).\)  Note that \(A_0 = \Id,\) so \(z_{ij}^{(0)}
\equiv \delta_{ij}.\) It is straightforward to check that the Hopf
algebra structure on the space of functions on $Maps_{1}(U,GL(2n,\C))$
(coming from the group structure on the gauge group), generated by the
$z_{ij}^{(M)}$, is the same as that on $Y_0(2n)$. (The identification
sends the element $z_{ij}^{(M)}$ to the generator
$t_{ij}^{(M)}$ of $Y_0(2n)$.) 

The Poisson structure on ${\cal G}_{2n}$ is given by the
first-order term in $h$ in the deformation of the algebra structure in
$Y_h(2n)$, so a glance at \eqref{eq:defYangcommrelations} yields the
Poisson structure 
\begin{equation}\label{eq:PBonG2n}
\{z_{ij}^{(M)}, z_{k\ell}^{(N)}\} = \sum_{r=0}^{\min(M,L)-1} 
(z_{kj}^{(r)} z_{i\ell}^{(M+L-1-r)} - 
z_{kj}^{(M+L-1-r)} z_{i\ell}^{(r)}).
\end{equation}
It remains to check that the product and Poisson structures on ${\cal
  G}_{2n}$ are compatible, i.e. the multiplication map \({\cal
  G}_{2n} \times {\cal G}_{2n} \to {\cal G}_{2n}\) is Poisson, where
  ${\cal G}_{2n} \times {\cal G}_{2n}$ has the product Poisson
  structure. This translates to the condition that for $f_1,f_2$
  functions on ${\cal G}_{2n}$ and $g, g' \in {\cal G}_{2n}$, 
\begin{equation}\label{eq:PLcompatibility}
\{f_1,f_2\}(gg') = \{L_g^* f_1, L_g^* f_2\}(g') + \{R_{g'}^* f_1,
R_{g'}^* f_2\}(g).
\end{equation}
This follows from the compatibility of the coalgebra and algebra
structures on $Y_h(2n)$ \eqref{eq:Hopfcompatibility}. In the case of
$Y_h(2n)$, we have defined $\Delta_1 \equiv 0$, so we get the
simplified compatibility equation
\begin{equation}\label{eq:Yangcompatibility}
\Delta(\mu_1(a_1 \otimes a_2)) = (\mu \otimes \mu_1 + \mu_1 \otimes
\mu)(\Delta^{13}(a_1)\Delta^{24}(a_2)).
\end{equation}
Here, $\mu$ indicates the commutative (pointwise) multiplication and
$\mu_1$, being the first-order term of the deformation of $\mu$, is the Poisson bracket. Finally,
$\Delta^{13}\Delta^{24}$ corresponds to the pullback induced by the
multiplication map
\begin{eqnarray*}
(m_{13}, m_{24}): &  G\times G\times G\times G  \rightarrow  G
  \times G \\
  &  (g_1, g_2, g_3, g_4) \mapsto (g_1g_3, g_2g_4)  \\
\end{eqnarray*}
From this it follows that the compatibility \eqref{eq:PLcompatibility}
is a consequence of \eqref{eq:Yangcompatibility}. 
\end{proof}

We now describe the classical limit of the twisted Yangian, which
turns out to be a Poisson homogeneous space associated to the
Poisson-Lie group ${\cal G}_{2n}$. We first define an involution
$\sigma$ on ${\cal G}_{2n}$ as follows. Let \(A(u) \in {\cal
  G}_{2n}.\) Then we define
\[
\sigma: A(u) \mapsto Q^{-1}(A(-u)^t)^{-1} Q,
\]
where $Q$ is the matrix defining the standard symplectic form ~\eqref{eq:symplform}. Note that this is just a point-wise version
of the standard involution on $GL(2n,\C)$ whose fixed point set is
$Sp(2n,\C)$. We then define \({\cal H}_{2n}\) to be the fixed point
set \({\cal G}_{2n}^{\sigma}\) under this involution.

\begin{theorem}\label{ph_ty}
The classical limit of the twisted Yangian is $\text{Fun}({\mathcal
G}_{2n}/{\mathcal H}_{2n})$.
\end{theorem}

\begin{proof}
We will show directly that the degeneration of the twisted Yangian, which
we denote by ${\cal A}$, consists exactly of functions $f$ on 
${\mathcal G}_{2n}$ with the property that for \(g,g' \in {\mathcal G}_{2n},\)
\begin{equation}\label{eq:GmodHrelation}
\quad g'=g \cdot h, \quad h \in {\mathcal H}_{2n} \Rightarrow f(g) = f(g').
\end{equation}
The generators of the twisted Yangian are given by 
the \(s_{ij}^{(M)}\) 
in~\eqref{eq:twY_gen}. This is
written collectively in matrix form as \(S(u) = T(u)T(-u)^{\tau}\) as in~\eqref{eq:Smatrixdef}. Here $\tau$ is the symplectic transpose defined in~\eqref{eq:sympltranspose}. 
Since the classical limit of $Y^{-}(2n)$ is generated by these
$s_{ij}^{(M)}$, it suffices to check the relation~\eqref{eq:GmodHrelation}
for these generators. 

Let $A(u) \in {\mathcal G}_{2n}$. We have 
\[
S(A(u)) := \left(\sum_M s_{ij}^{(M)}(A(u)) \cdot u^{-M} \right).
\]
It is straightforward to see that 
\[
S(A(u)) = A(u)A(-u)^{\tau}.
\]

We first show that the functions $S(u)$ are invariant under 
the action of ${\mathcal H}_{2n}$. Let 
\(B(u), C(u) \in {\mathcal G}_{2n},\) where \(B(u) = C(u)A(u)\) for an element
\(A(u) \in {\mathcal H}_{2n}.\) We want to show that 
\[
S(B(u)) = S(C(u)).
\]
By the above, this is equivalent to showing that
\[
B(u)B(-u)^{\tau} = C(u)C(-u)^{\tau}.
\]
Since \(B(u) =C(u)A(u)\) and the symplectic transpose $\tau$ 
is an antihomorphism, this is equivalent to showing that, for 
\(A(u) \in {\mathcal H}_{2n},\)
\[
A(u)A(-u)^{\tau} = \Id.
\]
By definition of ${\cal H}_{2n}$, 
we have 
\[
A(-u)^t = Q A(u)^{-1} Q^{-1}.
\]
Therefore
\begin{eqnarray*}
A(u)A(-u)^{\tau} &=& A(u) Q^{-1}A(-u)^t Q \\
 & = & A(u)Q^{-1}(Q A(u)^{-1} Q^{-1}) Q \\
 & = & A(u) A(u)^{-1} \\ 
 & = & \Id,
\end{eqnarray*}
and we are done. This argument is reversible, i.e. 
if \(B(u) = C(u)A(u),\) and \(S(B(u)) = S(C(u)),\) then \(A(u) 
\in {\mathcal H}_{2n}.\) Therefore the functions \(S(u) = (s_{ij}(u))\)
are precisely the functions on \({\mathcal G}_{2n}/{\mathcal H}_{2n}.\) 
\end{proof}

\begin{remark}
Heuristically, the fact that $Y^{-}(2n)$ has classical limit a
homogeneous space of ${\cal G}_{2n}$ may be motivated as follows. The
important observation is that $Y^{-}(2n)$, being a Hopf coideal of
$Y(2n)$, has a classical limit which is a Hopf coideal of
$Y_0(2n)$. Suppose now that $G$ is any Poisson-Lie group and $H$ is a subgroup.
The multiplication map \(m: G \times G \to G\) induces a map
\[
\overline{m}: G \times G/H \to G/H
\]
since $G/H$ is a quotient by $H$ on the right. This map then dualizes
to
\[
\text{Fun}(G/H) \to \text{Fun}(G) \otimes \text{Fun}(G/H),
\]
which is simply the coproduct \(\Delta: \text{Fun}(G) \rightarrow \text{Fun}(G) 
\otimes \text{Fun}(G)\) restricted to the subalgebra $\text{Fun}(G/H)$. 
Hence \(\Delta(\text{Fun}(G/H)) \in \text{Fun}(G) \otimes \text{Fun}(G/H),\)
and $\text{Fun}(G/H)$ is a Hopf coideal in $\text{Fun}(G)$. 
\end{remark}

Just as there is a geometric interpretation of the classical limit of
the Yangian, there is also a geometric interpretation of the classical
limit of the Yangian, namely, as a space of sections on
which ${\cal G}_{2n}$ acts transitively with stabilizer ${\cal
  H}_{2n}$. 
Let \(E = U \times \C^{2n}\) denote the total space of the trivial 
$\C^{2n}$-bundle over the formal neighborhood $U$. Let \(\alpha: U
\rightarrow U\) be the involution on $U$ given by
\[
\alpha: u \mapsto -u.
\]
Let \(\alpha^{*}E\) denote the pullback by $\alpha$ of the bundle $E$,
and let \((E \otimes \alpha^{*}E)^{*}\) be the bundle over $U$ whose
fiber over a point $u$ is the space of $\C$-bilinear pairings of
$E_{u}$ with \((\alpha^{*}E)_{u} = E_{-u}.\) We will show below that the
homogeneous space ${\mathcal G}_{2n}/{\mathcal H}_{2n}$ can be
identified with the subset of the space of sections $\Phi(u)$ of $(E
\otimes \alpha^{*}E)^{*}$ satisfying the ``skew'' condition
\begin{equation}\label{eq:skewCond}
\Phi(u)^t = - \Phi(-u)
\end{equation}
and the condition \(\Phi(\infty) = Q.\)
Here we have used the triviality of $E$ to identify all fibers with
a single $\C^{2n}$, and express \(\Phi(u) \in (E \otimes \alpha^{*}E)^{*}\) 
as a complex \(2n \times 2n\) matrix. 

\begin{theorem}\label{ty_geom} The classical limit of the twisted
  Yangian can be identified with sections of a vector bundle over the
  formal neighborhood $U$
  satisfying~\eqref{eq:skewCond} and \(\Phi(\infty) = Q.\) Thus
\[
{\mathcal G}_{2n}/{\mathcal H}_{2n} \cong \Gamma((E \otimes \alpha^{*}E)^{*})^{skew} :=
\{\Phi: \Phi(u)^t = - \Phi(-u), \Phi(\infty)=Q\} \subseteq 
\Gamma((E \otimes \alpha^{*}E)^{*}).
\]
\end{theorem}

\begin{proof}
We need to show that ${\mathcal G}_{2n}$ acts transitively on
\(\Gamma((E \otimes \alpha^{*}E)^{*})^{skew},\) where some point has 
stabilizer ${\mathcal H}_{2n}$. 
Let \(\Phi \in \Gamma((E \otimes \alpha^{*}E)^{*})\) be a 
bilinear pairing as above. The action of  
\(B \in {\mathcal G}_{2n}\) is given by 
\[
(B \cdot \Phi)(u) = (B(u)^{-1})^t \Phi(u) B(-u)^{-1}.
\]

We first show that there is an element of ${\mathcal G}_{2n}$ 
fixed by ${\mathcal H}_{2n}$. 
Consider the canonical section $\Omega$ of 
\((E \otimes \alpha^{*}E)^{*}\) given by the 
standard symplectic pairing $Q$ in~\eqref{eq:symplform}, so \(\Omega(u) \equiv Q.\) 
Now suppose \(A(u) \in {\mathcal G}_{2n},\) and 
\(A \cdot \Omega (u) = \Omega(u).\) This holds if and only
if 
\[
(A(u)^{-1})^t \Omega (u) A(-u)^{-1} = \Omega(u)
\]
as a formal power series in $u^{-1}$. 
This means 
\[
(A(u)^{-1})^t Q A(-u)^{-1} = Q.
\]
By properties of $Q$ and the definition of the symplectic transpose
$\tau$, the above is equivalent to the
equation
\[
A(u)A(-u)^{\tau} = \Id,
\]
i.e. \(A(u) \in {\mathcal H}_{2n}.\) Thus, the stabilizer of the ``constant'' pairing
\(\Omega(u) := Q\) is  precisely the subgroup ${\mathcal H}_{2n}$.

We now show that the action of ${\mathcal G}_{2n}$ is transitive. 
Let $\Phi$ be a bilinear pairing satisfying the conditions
\[
\Phi(u)^t = -\Phi(-u), \quad \Phi(\infty) = Q.
\]
We wish to show that there exists \(B(u) \in {\mathcal G}_{2n}\) such that 
\[
\Phi(u) = (B \cdot \Omega)(u) := (B(u)^{-1})^t Q B(-u)^{-1}.
\]
Equivalently, we wish to find \(C(u) = B(u)^{-1}\) such that
\[
\Phi(u) = C(u)^t Q C(u),
\]
as formal power series in $u^{-1}$. 

Let \(\Phi(u) = \sum_{k=0}^{\infty} \Phi_k u^{-k}\)
and set \(C(u) = \sum_{k=0}^{\infty} C_k u^{-k} \in {\mathcal G}_{2n}.\) 
By the conditions on $\Phi$, we must have that
\(\Phi_0 = Q,\) \(\Phi_k^t = \Phi_k\) for $k$ odd, and 
\(\Phi_k^t = -\Phi_k\) for $k$ even.
Expanding in powers of $u^{-1}$, we have
\[
C(u)^t Q C(-u) = \sum_{m=0}^{\infty} \left(\sum_{\substack{k+l=m\\k,l \geq 0}}
(-1)^l C_k^t Q C_l \right) u^{-m}.
\]
Thus, for each \(m \geq 0\) we wish to solve for $C_m$
in the equation
\[
\Phi_m = \sum_{\substack{k+\ell=m,\\ k,\ell\geq 0}} (-1)^\ell C_k^t Q C_\ell.
\]
For \(m=0,\) this is automatic since \(\Phi_0=Q, C_0 =\Id.\) The
equations for \(m \geq 1\) are solved for $C_m$ in a straightforward
manner by induction, and by using the fact that $\Phi_k$ is symmetric
for $k$ odd and antisymmetric for $k$ even.

\end{proof}

\subsection{The classical limit of the map $\Psi$}\label{subsec:CLmap}

We now turn our attention to the crucial algebra map $\Psi$ in
equation~\eqref{eq:MolevmapPsi}. In order to take the classical limit
of $\Psi$, it is convenient to first decompose $\Psi$ into pieces. 
We have the following: 
\begin{equation}\label{eq:defPsi}
\xymatrix{
Y^{-}(2) \ar @/^2pc/[rrr]^{\Psi} \ar[r]^{\imath} & Y^{-}(2n) \ar[r]^{\psi} & Y^{-}(2n)
\ar[r]^(.3){\phi} & U(\sp(2n,\C)) \\
}
\end{equation}
where \(\Psi := \phi \circ \psi \circ \imath.\) 
It turns out that
the image of $\Psi$ is in $\Ucent$ \cite{M4}.  The two maps $\imath$
and $\phi$ are relatively simple to describe. The technical heart of
this section lies in the analysis of the middle algebra map $\psi$.

We begin with the map $\imath$. The twisted Yangian $Y^{-}(2n)$ contains the
subalgebra $Y^{-}(2)$, generated by \(s_{i,j}(u), i, j \in \{-n,n\}.\)
We denote the inclusion map by 
\begin{equation}\label{twYinclusion}
\imath: Y^{-}(2) \into Y^{-}(2n).
\end{equation} 
We have chosen indices so that the subalgebra $Y^{-}(2)$ sit in the
``corner entries'' of the matrix $S(u)$.

In order to write down the maps $\phi$ and $\psi$, it is convenient to first
set some notation. Let ${\cal A}$ be an algebra with generators
$a_{ij}^{(M)}$. Let 
\[
A(u) = \sum_{M=0}^{\infty} A_M u^{-M}
\]
be a formal power series with matrix coefficients, where \(A_M := (a_{ij}^{(M)}).\) 
Then we may specify a map $f$ between ${\cal A}$ and an
algebra ${\cal B}$ by setting 
\[
f(A(u)) = B(u),
\]
where $f$ is understood to be linear over $u^{-1}$ and in the matrix
entries, so that 
\[
f(a_{ij}^{(M)}) = b_{ij}^{(M)},
\]
for some \(b_{ij}^{(M)} \in {\cal B}.\) We use this notation below. 

Let \(F = (F_{ij})\) be the \(2n \times 2n\) matrix whose $(i,j)$-th
entry is given by the generator $F_{ij}$
in $\sp(2n,\C)$~\eqref{sp2nbasis}. Then $\phi$
is defined by the formula
\begin{equation}\label{eq:psi}
\phi: S(u) \mapsto 1 + \frac{F}{u-\frac{1}{2}}.
\end{equation}
This is a well-defined algebra map \cite{MNO}.

We now come to the middle automorphism $\psi$ of $Y^{-}(2n)$. It turns out
that this is the restriction to $Y^{-}(2n)$ of an algebra
automorphism $\widehat{\psi}$ of $Y(2n)$ \cite{M4,MO},
so we have the commutative diagram
\[
\xymatrix{
Y(2n) \ar[r]^{\widehat{\psi}} & Y(2n)  \\
Y^{-}(2n) \ar[r]^{\psi} \ar @{^{(}->}[u] & Y^{-}(2n) \ar @{^{(}->}[u]\\
}
\]
and in order to take the classical limit of $\psi$, it is convenient to first
analyze that of $\widehat{\psi}$. 

We will describe $\widehat{\psi}$ as a composition of well-known ``basic''
algebra automorphisms of $Y(2n)$ \cite{MNO}. These are 
\smallskip
\begin{enumerate}
\item[1.] \(m_{g(u)}: T(u) \mapsto g(u)T(u)\), where $g(u)$ is a formal power 
series of the form 
\[
        g(u) := 1 + g_1 u^{-1} + g_2 u^{-2} + \ldots, \quad g_i \in \C.
\]
\item[2.] \(\tau_a: T(u) \mapsto T(u+a), a \in \C,\)
\item[3.] \(inv: T(u) \mapsto T(-u)^{-1},\)
\item[4.] \(\overline{\tau}: T(u) \mapsto T(-u)^{\tau},\)
\end{enumerate}
where $\tau$ is 
the symplectic transpose defined in~\eqref{eq:sympltranspose}. 

Using these ``basic'' automorphisms, the map $\widehat{\psi}$ is given
by 
\begin{equation}\label{eq:Y2naut}
\widehat{\psi} = \tau_n \circ \overline{\tau} \circ inv \circ 
m_{g(u)},
\end{equation}
where $g(u) = 1+g_1u^{-1}+g_2u^{-2}+\cdots$ is an element of
$Z(Y(2n))[[u^{-1}]]$, and $Z(Y(2n))$ is the center of $Y(2n)$ \cite{M4}. 

\begin{remark}
The proof in \cite{MNO} showing that $m_{g(u)}$ is an algebra
automorphism for $g(u)$ with coefficients in $\C$ also shows that
$m_{g(u)}$ is also an automorphism for any $g(u)$ with coefficients in
the center of the Yangian. This is because 
$g(u)T(u)$ still satisfies the ternary
relation~\eqref{eq:YangBaxter}. 
\end{remark}

We will now take the classical limits of these algebra maps $\phi,
\psi, \imath$. We first concentrate on the middle algebra map $\phi$. Since we
have decomposed $\phi$ as a composition of the basic
automorphisms, it suffices to calculate the classical limit of each of
these four basic types of algebra automorphisms. 

The theorem
below is the technical heart of this Section. The magic that
occurs here is that two of the basic automorphisms degenerate to
the {\em trivial} automorphism in the classical limit, while the other
two basic automorphisms remain the {\em same}. This
``all-or-nothing'' phenomenon accounts for the
amazing simplifications that occurs in the classical limit, and in
large part explains the simplicity of the formul{\ae} for the
functions integrating $\Ol$ in Theorem~\ref{thm:main}.

\begin{theorem}\label{thm:basicLimits}
The automorphisms \(m_{g(u)}, \tau_a\) degenerate to the identity 
at the level of the classical limit. The classical limits
of the automorphisms \(inv,
\overline{\tau}\)
are expressed by the same formul{\ae} as for the original
automorphisms. 
\end{theorem}

\begin{proof}
We will show that the multiplication map goes to the identity, and
that the inverse map remains the same. The calculations for the other
two types of automorphisms are similar.  

We first consider the multiplication map \(m_{g(u)}: T(u) \mapsto
g(u)T(u).\) In order to take the classical limit, we need to find an algebra
automorphism \((m_{g(u)})_h\) of $Y_h(2n)$ such that the following diagram commutes for
any \(h \neq 0\): 
\begin{equation}\label{eq:multgCommute}
\xymatrix{
Y_h(2n) \ar[d]_{\gamma_h} \ar[r]^{(m_{g(u)})_h} & Y_h(2n)
\ar[d]^{\gamma_h} \\
Y(2n) \ar[r]^{m_{g(u)}} & Y(2n) \\
}
\end{equation}
Here, $\gamma_h$ is the isomorphism between \(Y_h(2n)\) and \(Y(2n)\)
used in the proof of Lemma~\ref{lemma:deg_y}. We first observe that
the original automorphism \(m_{g(u)}\) acts as follows on the
generators:
\[
\xymatrix{
m_{g(u)}: \; t_{ij}^{(M)} 
\ar @{|->}[r] & {\displaystyle \sum_{\substack{K+L=M \\ K,L \geq 0}} g_K t_{ij}^{(L)}.
}
}
\]
Using this explicit form and the definition of $\gamma_h$ in the
proof of Lemma~\ref{lemma:deg_y}, one may immediately compute that the
map \(m_{g(u)} \circ \gamma_h\) takes 
\[
\xymatrix{
m_{g(u)}\circ \gamma_h:\;  \hs \tilde{t}_{ij}^{(M)} \ar @{|->}[r] & 
h^M \left({\displaystyle \sum_{\substack{K+L=M\\K,L\geq 0}}} g_K t_{ij}^{(L)}
\right). \\
}
\]
Thus, in order to have the diagram~\eqref{eq:multgCommute} commute, we
are forced to define \((m_{g(u)})_h\) as follows: 
\[
\xymatrix{
(m_{g(u)})_h:  \; \tilde{t}_{ij}^{(M)} \ar
  @{|->}[r] & 
{\displaystyle \sum_{\substack{K+L=M\\K,L\geq 0}}} g_K \tilde{t}_{ij}^{(L)} h^K.
}
\]
Thus when we set $h=0$, this degenerates to the map
\[
\xymatrix{
(m_{g(u)})_{0}: \hs \tilde{t}_{ij}^{(M)} \ar @{|->}[r] &
\tilde{t}_{ij}^{(M)},
}
\]
i.e. the identity map. 

Now we consider the automorphism \(inv: T(u) \mapsto T(-u)^{-1}.\) We
wish to make a diagram similar to~\eqref{eq:multgCommute} commute,
using now $inv$ instead of $m_{g(u)}$. We first observe that the
automorphism $inv$ behaves as follows on the generators, \(M \geq 1\):
\[
\xymatrix{
inv: \; t_{ij}^{(M)} \ar @{|->}[r] & {\displaystyle \sum_{s=1}^M (-1)^{M+s}} \left({\displaystyle \sum_{\substack{m_1+\ldots +m_s=M,\\
m_i \geq 1}} \sum_{\substack{a_1, \ldots, a_{s-1}\\ a_i \in {\cal I}}} t_{i,a_1}^{(m_1)} t_{a_1, a_2}^{(m_2)} 
\cdots t_{a_{s-1}, j}^{(m_s)}} \right).
}
\]
The generators with \(M=0\) are sent to themselves. 
Then we immediately compute that for the map 
\(inv \circ \gamma_h,\) we have for \(M \geq 1\) 
\[
\xymatrix{
inv \circ \gamma_h: \hs \tilde{t}_{ij}^{(M)}\ar @{|->}[r] &
 h^M \left[ {\displaystyle \sum_{s=1}^M (-1)^{M+s}}
\left( {\displaystyle 
\sum_{\substack{m_1 + \cdots + m_s =M \\ m_i \geq 1}} 
\sum_{\substack{a_1, \ldots, a_{s-1}\\ a_i \in {\cal I}}}}
t_{i,a_1}^{(m_1)} \ldots t_{a_{s-1},j}^{(m_s)} \right)  \right].
}
\]
Thus we see that the map \(inv_h\) on $Y_h(2n)$ is 
{\em independent} 
of $h$, and thus the classical limit $inv_0$ is given by the same formula as for
the map $inv$ on $Y(2n)$.
\end{proof}

The computation of the classical limits of the basic automorphisms
immediately leads us to a quick computation of the classical limit of both
$\widehat{\psi}$ and $\psi$. 

\begin{corollary}\label{cor:hatPsiCL}
The map \(\widehat{\psi}\) 
degenerates to \(\overline{\tau} \circ inv\) in 
the classical limit, and sends 
\[T(u) \mapsto (T(u)^{\tau})^{-1}.\] 
\end{corollary}

\begin{proof}
The map \(\widehat{\psi}\) is decomposed as \(\tau_n \circ \overline{\tau}
\circ inv \circ m_{g(u)}.\) Both
$m_{g(u)}$ and $\tau_n$ degenerate to the identity,
and $inv$ and $\overline{\tau}$ remain the same. Thus the limit is simply
$\overline{\tau} \circ inv$, as desired. 
\end{proof}

\begin{remark}
On the quantum level, the map \(T(u) \mapsto (T(u)^{\tau})^{-1}\) is
also a coalgebra automorphism \cite{MNO}. On the classical level, one
therefore expects it to degenerate to a group automorphism of ${\cal
  G}_{2n}$. This is indeed the case, as may be checked directly. 
\end{remark}

The classical limit $\widehat{\psi}_0$ of $\widehat{\psi}$ preserves the subgroup ${\cal
  H}_{2n}$, and therefore induces a map on the homogeneous space
  \({\cal G}_{2n}/{\cal H}_{2n},\) giving us the classical limit of
  $\psi$. We now give the explicit formula for this map, written in
  terms of the coordinates \(S(u) = (s_{ij}(u)).\) 

\begin{corollary}\label{cor:psiCL}
The classical limit $\psi_0$ of the map $\psi$ is given by 
\(S(u) \mapsto (S(-u))^{-1}.\) 
\end{corollary}

\begin{proof}
Since $\widehat{\psi}_0$ takes \(T(u) \mapsto (T(u)^{\tau})^{-1},\) 
we see that
\(T(-u)^{\tau} \mapsto ((T(-u)^{\tau})^{-1})^{\tau}.\) 
Hence the coordinates \(S(u) = 
T(u)T(-u)^{\tau}\) are mapped to \((T(u)^{\tau})^{-1}
((T(-u)^{\tau})^{-1})^{\tau}.\) Since $\tau$ is an involutory 
automorphism, we have that \((T(-u)^{\tau})^{-1} = (T(-u)^{-1})^{\tau},\)
and thus 
\[
S(u)=T(u)T(-u)^{\tau} \mapsto (T(u)^{\tau})^{-1} T(-u)^{-1} = S(-u)^{-1}.
\]
Thus the map on the homogeneous space is given by 
\(S(u) \mapsto S(-u)^{-1}.\) 
\end{proof}

Now that we have calculated the limit $\psi_0$ of $\psi$, we now focus on the two algebra maps $\phi$ and $\imath$
in~\eqref{eq:defPsi}. As advertised previously, the calculations here
are more straightforward than those for $\psi$. 

\begin{proposition}\label{prop:phiCL}
The classical limit $\phi_0$ of the map $\phi$ is given by
\[
\phi_0: 
\begin{cases}
\tilde{s}_{ij}^{(0)} := \delta_{ij} \mapsto  \delta_{ij}, \\
\tilde{s}_{ij}^{(1)} \mapsto F_{ij}, \\
\tilde{s}_{ij}^{(M)} \mapsto 0, \quad \mbox{for} \hs M \geq 2. \\
\end{cases}
\]
\end{proposition}

\begin{proof}
In order to compute
$\phi_0$, we need to find algebra maps $\phi_h$ such that a diagram similar
to~\eqref{eq:multgCommute} commutes. Using the expansion
\[
\frac{1}{u - \frac{1}{2}} = \sum_{r=0}^{\infty} (-1)^r \left(\frac{1}{2}\right)^r
u^{-r-1},
\]
we see that an explicit formula for $\phi$ is given by 
\[
 s_{ij}^{(M)} (-1)^{M-1} \left(\frac{1}{2}\right)^{M-1} \mapsto F_{ij}.
\]
It is then a straightforward computation to see that $\phi_h$ is given
by
\[
  \tilde{s}_{ij}^{(M)} \mapsto F_{ij} \left(\frac{1}{2}\right)^{M-1} (-1)^{M-1}
h^{M-1}.
\] 
In particular, for \(h=0,\) the only non-zero images are those of
$\tilde{s}_{ij}^{(1)}$, and we have the formul{\ae} as desired. 
\end{proof}

In order to give a geometric description of $\phi_0$, we first recall
that the classical limit of  \(U_h(\sp(2n,\C))\) is well-known to be 
\(Sym(\sp(2n,\C)^{*})\) (see e.g. \cite{CP}), the space of 
polynomial functions on the dual of the Lie algebra $\sp(2n,\C)$. 
The underlying Poisson space is \(\sp(2n,\C)^{*}\).  With this in
hand, it is immediate from Proposition~\ref{prop:phiCL} that 
the geometric map corresponding to $\phi_0$ is given by 
\begin{equation}\label{eq:phi_cl}
        X \in \sp(2n,\C)^{*} \mapsto \Id + X u^{-1},
\end{equation}
where the image is interpreted as representing an element in ${\cal G}_{2n}/{\cal
  H}_{2n}$. 

Finally, we calculate the classical limit of the inclusion map
\(\imath\). 

\begin{proposition}\label{prop:imathCL}
The classical limit $\imath_0$ of the map $\imath$ is given by 
\[
\tilde{s}_{\pm n, \pm n}^{(M)} \in Y^{-}(2) \longrightarrow
\tilde{s}_{\pm n, \pm n}^{(M)} \in Y^{-}(2n).
\]

\end{proposition}

\begin{proof}
The classical limit $\imath_0$ is given by the same formula as for $\imath$ since
the inclusion map preserves degrees (where
\(\mathrm{deg}(s_{ij}^{(M)}) = M\)). 
\end{proof}

Interpreted as a geometrical map from \({\cal G}_{2n}/{\cal H}_{2n}\)
to \({\cal G}_2/{\cal H}_2,\) it is given in coordinates by 
\[
S(u) \mapsto \left[ \begin{array}{cc} s_{-n,-n}(u) & s_{-n,n}(u) \\
    s_{n,-n}(u) & s_{n,n}(u) \\ \end{array} \right], 
\]
i.e. it is the quotient map given by ``taking the corner entries'' of
$S(u)$. For an element
\(A(u) \in {\cal G}_{2n}/{\cal H}_{2n},\) let \(A(u)_{\pm n,\pm n}\) denote the
element in ${\cal G}_2/{\cal H}_2$ gotten by taking the corner entries
as above.

We now give the formula for the classical limit of $\Psi$, interpreted
as a geometric map on the underlying Poisson spaces. Recall that the
Poisson space which is the classical limit of $U(\sp(2n,\C))$ is
$\sp(2n,\C)^*$ \cite{CP}, and thus the classical limit of the
centralizer \(\Ucent\) is the Poisson quotient of the Lie algebra dual
by a subgroup,
\(\sp(2n,\C)^*/Sp(2(n-1),\C).\) 

\begin{theorem}\label{thm:PsiCL}
The classical limit $\Psi_0$ of $\Psi$ is given by the composition $\imath_0
\circ \psi_0 \circ \phi_0$. As a map on the underlying geometric
objects, $\Psi_0$ is given in coordinates as follows. Let \(X \in
\sp(2n,\C)^*.\) 
\[
\Psi_0:  \hs X \mapsto \left(\sum_{M=0}^{\infty} X^M
u^{-M}\right)_{\pm n, \pm n}
\]
This is invariant under the action of $Sp(2(n-1),\C)$, so it is
well-defined on $\sp(2n,\C)^*/Sp(2(n-1),\C)$. 
\end{theorem}

\begin{proof}
The first map $\phi_0$ sends $X$ to \(\Id + Xu^{-1}.\) The map
$\psi_0$ sends
\[
\Id + Xu^{-1} \mapsto (\Id - Xu^{-1})^{-1} = \Id + Xu^{-1} + X^2
u^{-2} + \cdots,
\]
and then $\imath_0$ takes the corner entries, as desired. Since the
\((\pm n, \pm n)\) entries are untouched by the action of
$Sp(2(n-1),\C)$, the map is invariant under this action. 
\end{proof}

\subsection{The derivation of the integrable system}\label{subsec:Derivation}

So far, in taking the geometric classical limit of Molev's
constructions in \cite{M4}, we have followed his conventions in using
algebras over $\C$. In so doing, we have obtained, in
Theorem~\ref{thm:PsiCL}, a map on the quotient of the complex Lie
algebra dual $\sp(2n,\C)^*$. In this section, we take the compact
analogue of $\Psi_0$ in order to obtain the symplectic geometric
picture, with $\R$-valued functions on (a quotient of) the compact
form $\unH^*$. 

Recall that \(\gl(n,\H)\) is naturally a subalgebra of $\gl(2n,\C)$ by
a restriction of scalars. Given an element \(A(u) \in
\gl(n,\H)[[u^{-1}]],\) we denote by \(A(u)_{nn}\) the element of
\(\gl(1,\H)[[u^{-1}]] \cong \H[[u^{-1}]]\) obtained by taking the $(n,n)$-th matrix entry of
$A(u)$.

\begin{definition}\label{def:PsiH}
For \(X \in \unH \cong \unH^{*},\) we define the {\em compact form $\Psi_{\H}$}
of $\Psi_0$ as 
\begin{equation}\label{eq:PsiH} 
\Psi_{\H}: X \mapsto (1 + X u^{-1} + X^2 u^{-2} + \ldots)_{nn}.
\end{equation}
This map is well-defined on the quotient $\unH^{*}/U(n-1,\H)$ since it is 
invariant under the action of $U(n-1,\H)$. 
\end{definition}

As we saw in
Section~\ref{subsec:GCforUnC}, the key ingredient for the construction
of a Gel$'$fand-Cetlin basis for the case of $U(n,\C)$ is the presence of a large family of
commuting operators on $V(\lambda)$, the classical limit of which gave
a large family of Poisson-commuting functions on $\Ol$. Following this
example, we look now for commuting elements of $Y^{-}(2)$, which via
Molev's map $\Psi$ are commuting operators on an irreducible
representation $V(\lambda)$ of $U(n,\H)$. We are in luck: Molev
observes \cite{M4} that the coefficients of \(\mathrm{tr}(S(u))\)
generate a commuting subalgebra of $Y^{-}(2)$. Thus, in the classical
limit, the functions \(s_{-n,-n}^{(M)} + s_{n,n}^{(M)}\) for \(M \geq
1\) Poisson-commute on \({\cal G}_2/{\cal H}_2,\) and hence (since
$\Psi_{\H}$ is a Poisson map) their pullbacks Poisson-commute on
$\unH^*/U(n-1,\H)$. 
Note that as functions on $\sptyt$, the
$s_{-n,-n}^{(M)} + s_{n,n}^{(M)}$ are a priori $\C$-valued, since they
simply read off certain matrix entries in $\gl(2n,\C)$. However,
the image of $\Psi_{\H}$ is by definition contained in the
intersection 
\(\gl(1,\H)[[u^{-1}]] \cap {\cal G}_2/{\cal H}_2 \cong \H[[u^{-1}]]
\cap {\cal G}_2/{\cal H}_2.\) 
Restricted to this subset, the functions are in fact $\R$-valued. We
record the following calculation.

\begin{lemma}\label{lemma:realValued}
The functions \(\bigg(s_{-n,-n}^{(M)} + 
s_{n,n}^{(M)}\bigg)\) are $\R$-valued when restricted to
\(\H[[u^{-1}]] \cap {\cal G}_2/{\cal H}_2.\)  
In particular, for \(A(u) \in {\cal G}_2/{\cal H}_2,\)
\[
\bigg(s_{-n,-n}^{(M)} + s_{n,n}^{(M)} \bigg)(A(u)) = 2 \cdot \mathrm{Re}(A_M),
\]
where \(A_M \in \H \subset \gl(2,\C),\) and $\mathrm{Re}(A_M)$ 
denotes the real part in the quaternionic sense. 
\end{lemma}

\begin{proof}
An element in \(\H[[u^{-1}]] \cap {\cal G}_2/{\cal H}_2\) is a
formal series
\[
A(u) = {\Id} + A_1 u^{-1} + A_2 u^{-2} + \ldots \in {\Id} + 
u^{-1} \gl(2,\C)[[u^{-1}]],
\]
where \(A(u) = B(u)B(-u)^{\tau},\) for \(B(u) \in {\cal G}_2.\) Since
\(A(u) \in \H[[u^{-1}]],\) each coefficient \(A_M\) is an element in $\H$, where
these are considered as element of $\gl(2,\C)$ by the standard
inclusion \(\H \into \gl(2,\C)\):
\[
A_M = \alpha_M + j\beta_M \mapsto 
\left[
\begin{array}{cc} \overline{\alpha}_M & \beta_M \\
                -\overline{\beta}_M & \alpha_M \\
\end{array} 
\right].
\]
Then it is immediate that 
\[
\bigg(s_{-1,-1}^{(M)}+s_{1,1}^{(M)}\bigg)(A(u)) = \overline{\alpha}_M + \alpha_M = 
2 \cdot \mathrm{Re}(\alpha),
\]
and in particular is $\R$-valued. Here, $\mathrm{Re}(\alpha)$ denotes
the real part in the $\H$ sense. 
\end{proof}

We are now prepared to obtain the formul{\ae} for the functions $f_{n,m}$
used in Section~\ref{subsec:main}. 

\begin{theorem}\label{thm:fmnDerive}
Let \(\Ol\)
be a generic coadjoint orbit of $U(n,\H)$.  
Let \(\Psi_{\H}: \unH/U(n-1,\H) \to {\cal G}_2/{\cal H}_2\) denote the compact
form of \(\Psi_0.\) 
Let \(X \in \Ol,\) so \(X = A D_{\lambda} A^{*}\) for some \(A \in U(n,\H).\) 
Then 
\begin{equation}\label{eq:fmnDerive}
(\Psi_{\H})^{*} \bigg(s_{-1,-1}^{(M)} + s_{1,1}^{(M)} \bigg)(X) =
 \mathrm{rtr}((AD_{\lambda}A^*)^M E_{nn}), 
\end{equation}
where the $a_{n,\ell}$ denote the entries in the bottom row of 
\(A = (a_{ij}).\) 
In particular, for $M$ odd, the pullback functions are identically
0. For all $M$, the pullback functions are invariant under $U(n-1,\H)$
and hence well-defined on the quotient \(\unH^{*}/U(n-1,\H).\) 
\end{theorem}

\begin{proof}
Given the diagonalization \(X = A D_{\lambda} A^*,\) any $X^M$ is of
the form \(X^M = (AD_{\lambda}A^*)^M = A D_{\lambda}^M A^*.\) Since
$D_{\lambda}$ is diagonal, its power \(D_{\lambda}^{M}\) is 
a diagonal matrix with diagonals given by 
\[
((i\lambda_1)^M, \ldots, 
(i\lambda_n)^M).
\]
The $(n,n)$-th entry in $X^M$ is therefore
given by 
\[
(X^M)_{(n,n)} = \sum_{\ell=1}^n a_{n,\ell} (i\lambda_\ell)^m \overline{a_{n,\ell}},
\]
where \(A = (a_{ij}).\) Since the functions
$s_{-1,-1}^{(M)}+s_{1,1}^{(M)}$ read off the real part 
of $M$-th coefficient, as in Lemma~\ref{lemma:realValued},
\[
(\Psi_{\H})^* \bigg(s_{-1,-1}^{(M)}+s_{1,1}^{(M)} \bigg)(X) = 2 \cdot \sum_{\ell=1}^n 
\mathrm{Re} \bigg(a_{n,\ell} (i\lambda_\ell)^M \overline{a_{n,\ell}} \bigg).
\]
Using the reduced trace~\eqref{redtr}, this can be rewritten as  
\begin{equation}\label{eq:fmnReducedTrace}
(\Psi_{\H})^{*} \bigg(s_{-1,-1}^{(M)} + s_{1,1}^{(M)} \bigg)(X) =
  \mathrm{rtr}((AD_{\lambda}A^*)^M E_{nn}).
\end{equation}
Notice that for $M$ odd, each $(i\lambda_l)^M$ is pure imaginary, and
hence each term in the sum is pure imaginary in $\H$.
Hence its real part is $0$. Therefore, for
$M$ odd, this pullback is identically 0 as a function on $\Ol$. Hence
we only get non-trivial functions for $M=2m$ even.  Moreover, since
the function reads off only the $(n,n)$-th entry, it is invariant
under $U(n-1,\H)$. 
\end{proof}

\appendix

\renewcommand{\theequation}{A.\arabic{equation}}
\setcounter{equation}{0}

\section*{Appendix A: On quaternionic linear algebra and $U(n,\H)$}\label{app:UnH}
\addcontentsline{toc}{section}{Appendix A: On quaternionic linear algebra and $U(n,\H)$}

The quaternions $\H$ are defined as the set of quadruples
\[
q = a + ib + jc + kd,
\]
where \(a,b, c, d \in \R\) and the \(i,j,k\) satisfy the 
relations \(i^2 = j^2 = k^2 = -1, ij=k.\)
The quaternions are {\em not}
commutative, since \(ij = -ji = k.\)
We define conjugation in $\H$ as 
\[
\overline{q} = a - ib - jc - kd,
\]
for $q$ as above. We define
\(\text{Re}(q):= a\) and \(\text{Im}(q) := ib+ jc + kd\) 
for $q$ as above. 
We define the norm of an element $q$ to be 
\(\|q\| = \sqrt{q \overline{q}} \in \R.\)

We now describe our conventions for linear algebra over $\H$. 
Let $\H^n$ be the $n$-dimensional quaternionic vector space of 
$n$-tuples in $\H$, equipped with scalar multiplication
by $\H$ on the {\em right}. 
Elements of $\H^n$ will be represented by column vectors, and
$\H$-linear transformations will then be represented by matrix
multiplication on the {\em left}. We denote by $\gl(n,\H)$ the algebra
of \(n \times n\) matrices with entries in $\H$.
The standard
basis vectors are the \(e_i = (0,0, \cdots, 0,1,0,\cdots,0)^t,\) where
the $1$ is in the $i$-th place. We define the conjugate $\overline{v}$ of a
vector $v \in \Hn$ componentwise, and the norm $\|v\|$ also as usual,
by a sum of norms of the components. For any \(m \times n\) matrix $A
=(a_{ij})$ with entries in $\H$, we define the conjugate transpose as
the $n \times m$ matrix
\[
A^{*} := (\overline{A})^t,
\]
where conjugation is quaternionic conjugation, and the
transpose $t$ is as usual. 

Given two vectors \(v = (v_1, \ldots, v_n)^t, w= (w_1, \ldots, w_n)^t
\in \H^n,\) the standard quaternionic hermitian form 
is defined by 
\[
\langle v,w \rangle := \sum_{l=1}^n \overline{v}_l w_l.
\]
More compactly, 
\[
\langle v,w\rangle := v^{*}w,
\]
where the conjugate transpose is defined above.
We define the {\em compact symplectic group} to be the subset of
$\H$-linear transformations preserving the quaternionic hermitian
form. 
\[
U(n,\H) := \{ A \in \gl(n,\H): \langle Av, Aw\rangle = \langle v,w\rangle, \forall v,w \in \H^n\}.
\]
Again, more compactly, 
\[
U(n,\H) = \{A \in \gl(n,\H): A^{*} A = \Id \}.
\]
The Lie algebra is 
\[
\unH = \{X \in \gl(n,\H): X^{*} + X = 0\}.
\]

Observe that \(\H \cong \C \oplus j\C,\) where the $\C$ is thought of
as \(\R \oplus i\R.\) Similarly, \(\H^n \cong \C^{n} \oplus j\C^n.\) 
By a restriction of scalars from $\H$ to $\C$, $\Hn$ can be thought
of as a $\C$-vector space of dimension $2n$, with 
ordered basis 
\begin{equation}\label{basis}
\{e_{-n}, e_{-(n-1)}, \ldots, e_{-2},e_{-1}, e_1, e_2, \ldots,
e_{n-1},e_n\},
\end{equation}
where \(e_{-i} := e_i \cdot j \in \Hn \cong \C^{2n}.\) With this basis
of $\C^{2n}$ in mind, we take the indexing set to be
\begin{equation}\label{eq:indexingSet}
{\cal I} := \{-n, -n+1, \ldots, -2, -1, 1, 2, \ldots, n-1, n\}.
\end{equation}
Note that the index $0$ is skipped. 
Since any $\H$-linear map is also $\C$-linear, there is a natural 
inclusion $\gl(n,\H)$ into $\gl(2n,\C)$, the space of $2n \times 2n$
matrices with $\C$ entries. 

Using the decomposition \(\H \cong \C \oplus j \C,\) we may also write
the quaternionic hermitian form as a sum
\[
\langle v,w\rangle = H(v,w) + j Q(v,w) \in \C \oplus j\C,
\]
where $H$ is the standard hermitian form on $\C^{2n}$, and $Q$ is
a {\em complex symplectic form} on $\C^{2n}$. Written with respect to
the ordered basis of $\C^{2n}$ above, the symplectic form is given by
\[
Q(v,w) := v^t Q w.
\]
Here, $Q$ is the \(2n \times 2n\) matrix
\begin{equation}\label{eq:symplform}
Q = \left[ \begin{array}{cc} \mathbf{0} & \tilde{I}_n \\
           -\tilde{I}_n & \mathbf{0} \end{array} \right],
\end{equation}
where \(\tilde{I}_n\) is the \(n \times n\) matrix with
ones along the {\em antidiagonal}, i.e. \(\tilde{I}_n = (a_{ij})\) where
\(a_{ij} = \delta_{i, (n+1)-j}.\)

Using the symplectic form $Q$ above, we define the {\em symplectic
  transpose} $\tau$ as the involution on $\gl(2n,\C)$ which satisfies 
\[
Q(Av,w) = Q(v,A^{\tau}w)
\]
for all \(v,w \in \C^{2n}, A \in \gl(2n,\C).\) An explicit formula for
the symplectic transpose is given by 
\begin{equation}\label{eq:sympltranspose}
\tau: A \mapsto Q^{-1} A^t Q,
\end{equation}
where the matrix $Q$ is given in \ref{eq:symplform} and the $t$ is the
usual transpose. 

The subgroup
of $\gl(2n,\C)$ preserving the complex symplectic form $Q$ is defined to
be the {\em complex symplectic group} $Sp(2n,\C)$. We have
\begin{equation}\label{def1_Sp}
Sp(2n,\C) := \{ A \in \gl(2n,\C): Q(Av,Aw) = Q(v,w) \hs \hs \forall v,w \in
\C^{2n}\},
\end{equation}
which again can be rewritten as 
\begin{equation}\label{def2_Sp}
Sp(2n,\C) = \{A \in \gl(2n,\C): A^t Q A = Q \}
\end{equation}
for the matrix $Q$ defined above. Since the compact symplectic group
is precisely the subgroup preserving both the usual hermitian form on
$\C^{2n}$ the symplectic form $Q$, we have
\[
U(n,\H) = U(2n,\C) \cap Sp(2n,\C).
\]
We may conclude from this fact, plus a dimension count, 
that $U(n,\H)$ is the compact form of
$Sp(2n,\C)$. This justifies the terminology.

The Lie algebra of $Sp(2n,\C)$ is given by 
\[
\sp(2n,\C) = \{ X \in \gl(2n,\C): X^t Q + QX = 0 \}.
\]
A basis of $\sp(2n,\C)$ is given by the elements 
\begin{equation}\label{sp2nbasis}
        F_{i,j} = E_{i,j} + \mathrm{sgn}(i) \cdot \mathrm{sgn}(j) E_{-j, -i},
\end{equation}
where the index set for basis elements of $\gl(2n,\C)$ is
${\cal I} := \{ -n, -(n-1), \ldots, -1, 1, 2, \ldots, n\}$ (note that we skip the
index $0$). 

We now collect some standard Lie-group-theoretic facts
about the compact symplectic group $U(n,\H)$. A further discussion can
be found in \cite{Cur}. 

As a real manifold, $U(n,\H)$ has dimension $2n^2 + n$. 
The dimension of
a maximal torus of $U(n,\H)$ is $n$. We will always take as choice of maximal
torus the diagonal subgroup $T^n \cong (S^1)^n$, where \(S^1 =
\{e^{i\theta}\}.\) This is the $S^1$ sitting in the {\em first} factor
of \(\H \cong \C \oplus j \C.\) 
Any element $X$ of the Lie algebra $\unH$
can be conjugated by an element $A$ of $U(n,\H)$ to a diagonal matrix:
\[
A X A^{*}  =  D_{\lambda}. 
\]
Since the Weyl group of $U(n,\H)$ is the group of {\em signed}
permutations, we may choose $D_{\lambda}$ to be of the form 
\begin{equation}\label{eq:UnHdiagonal}
D_{\lambda} := 
\left[
\begin{array}{cccc}
i\lambda_1 & & & \\
  & i \lambda_2 & & \\
  &   & \ddots & \\
  &   &    & i\lambda_n \\
\end{array}
\right], \qquad \mbox{where} \quad 0 \geq \lambda_1 \geq \lambda_2
\geq \cdots \geq \lambda_n.
\end{equation}
We denote by \(\lambda = (\lambda_1, \ldots, \lambda_n)\) the
$n$-tuple of ``eigenvalues.'' 

The Killing form on $\unH$ is the restriction of the usual Killing
form on $\u(2n,\C)$ restricted to $\unH$.
This is called the {\em reduced trace} pairing, where the reduced trace of a
quaternionic matrix $A$ is defined by
\begin{equation}\label{reducedtrace}
\text{rtr}(A) := tr(\iota(A)),
\end{equation}
and \(\iota\) is the inclusion \(\unH \into \u(2n,\C).\) 
It will be convenient to express the 
reduced trace purely in quaternionic terms.
For an element 
\(A = (a_{ij}) \in \gl(n,\H),\) the reduced trace is given by
\begin{equation}\label{redtr}
\text{rtr}(A)  =  2 \cdot \text{Re}\big(\sum_{i=1}^n a_{ii}\big),
\end{equation}
i.e. it is twice the real part of the sum of the diagonals. 

\renewcommand{\theequation}{B.\arabic{equation}}
\setcounter{equation}{0}

\renewcommand{\thetheorem}{B.\arabic{theorem}}
\setcounter{theorem}{0}

\section*{Appendix B: The Yangian and twisted Yangian}\label{app:Yangian}
\addcontentsline{toc}{section}{Appendix B: The Yangian and twisted Yangian}

We now briefly recall essential facts about the 
Yangian and twisted Yangian. See \cite{MNO} for
details. 

The Yangian \(Y(2n) = Y(\gl(2n))\) is defined\footnote{The Yangian
  $Y(N)$ can be defined for any $N$, but we are only interested in the
  case $N=2n$ even.} as a complex associative unital algebra with
  countably many generators \(t_{ij}^{(1)}, t_{ij}^{(2)}, \ldots,\)
  where the indices \(i, j\) are in the indexing set ${\cal I}$ as in
  Appendix A. In particular, the Yangian is
  infinite-dimensional.  The generators satisfy the following defining
  commutation relations
\begin{equation}\label{Yangianrelations}
[t_{ij}^{(M)}, t_{kl}^{(L)}] = \sum_{r=0}^{\min(M,L)-1}
\left( t_{kj}^{(r)} t_{il}^{(M+L-1-r)} - t_{kj}^{(M+L-1-r)}t_{il}^{(r)} 
\right).
\end{equation}
Moreover, the Yangian is a Hopf algebra with coproduct defined as 
\[
\Delta(t_{ij}(u)) := \sum_{a=-n}^n t_{ia}(u) \otimes t_{aj}(u).
\]

It turns out to be useful to write these relations more compactly. 
We first set some definitions and notation. Define 
the formal power series 
\[
t_{ij}(u) := \delta_{ij} + t_{ij}^{(1)}u^{-1} + t_{ij}^{(2)} u^{-2} + 
\cdots \in Y(2n)[[u^{-1}]],
\]
and assemble them in a single ``$T$-matrix'' as
\[
T(u) := \sum_{i,j} t_{ij}(u) \otimes E_{ij} \quad \in \quad 
Y(2n)[[u^{-1}]] \otimes
\mathrm{End}(\C^{2n}).
\]
Moreover, for an operator \(X \in \mathrm{End}(\C^{2n}),\) we set 
\[
X_1 := X \otimes \Id, \quad X_2 := \Id \otimes X,
\]
in \(\mathrm{End}(\C^2n \otimes \C^{2n}).\) Then we may define 
\[
T_1(u) := \sum_{i,j}^{2n} t_{ij}(u) \otimes (E_{ij})_1, \quad T_2(v)
:= \sum_{i,j}^{2n} t_{ij}(v) \otimes (E_{ij})_2.
\]
In a similar spirit, we define
\[
T_{[1]}(u) := \sum_{i,j=1}^{2n} t_{ij}(u)
\otimes 1 \otimes E_{ij}
\in Y(2n)[[u^{-1}]]^{\otimes 2} \otimes \mathrm{End}(\C^{2n}),
\]
and 
\[
T_{[2]}(u) := \sum_{i,j=1}^{2n} 1 \otimes t_{ij}(u) \otimes E_{ij}
\in Y(2n)[[u^{-1}]]^{\otimes 2} \otimes \mathrm{End}(\C^{2n}).
\]

With the notation above, the relations in equation~\eqref{Yangianrelations} 
can be written more compactly as the following single
relation for the $T$-matrix: 
\begin{equation}\label{eq:YangBaxter}
R(u-v) T_1(u)T_2(v) = T_2(v)T_1(u)R(u-v),
\end{equation}
where \(R(u) := 1 - \frac{P}{u},\) and $P$ is the permutation
operator on \(\C^{2n} \otimes \C^{2n}.\) \footnote{The object \(R(u) \in 
End(\C^{2n}) \otimes End(\C^{2n}) \otimes \C(u)\) is called the {\em 
Yang-Baxter $R$-matrix.} The relation above is called the {\em ternary 
relation.}} The coproduct structure may similarly be described by the
single matrix equation
\begin{equation}\label{coalgdefTmatrix}
\Delta(T(u))= T_{[1]}(u) T_{[2]}(u) \in Y(2n)[[u^{-1}]]^{\otimes 2}
\otimes \mathrm{End}(\C^{2n}).
\end{equation}

The antipode map is defined by \(S: T(u) \mapsto T^{-1}(u)\) and the
counit is \(\epsilon(T(u)) := 1.\) This makes the Yangian into a Hopf
algebra \cite{MNO}. 

The twisted Yangian $Y^{-}(2n)$ is defined as the subalgebra of $Y(2n)$ 
with generators the entries in the matrix 
\begin{equation}\label{eq:Smatrixdef}
S(u) := T(u)T(-u)^{\tau},
\end{equation}
where $\tau$ was defined in~\eqref{eq:sympltranspose}.
In terms of the matrix entries
\[
s_{ij}(u) = \sum_{a=-n}^n \theta_{aj} t_{ia}(u)t_{-j,-a}(-u) =
\sum_{M=0}^{\infty} s_{ij}^{(M)}u^{-M}, 
\]
we obtain the formula for the generators $s_{ij}^{(M)})$ in terms of
the generators of $Y(2n)$. 
\begin{equation}\label{eq:twY_gen} 
s_{ij}^{(M)}  =  \sum_{a=-n}^n \sum_{K+L=M, K,L \geq 0} (-1)^L \theta_{aj}
t_{ia}^{(K)} t_{-j,-a}^{(L)}.
\end{equation}

It can be shown (\cite{MNO}, Prop 4.17) that the twisted Yangian is a left
Hopf coideal of the Yangian, i.e.
\[
\Delta(Y^{-}(2n)) \subset Y(2n) \otimes Y^{-}(2n).
\]
In particular, 
\[
\Delta(s_{ij}(u)) = \sum_{k,l} \theta_{lj} t_{ik}(u) t_{-j,-l}(-u) 
\otimes s_{kl}(u),
\]
where \(i,j \in \{-n,-(n-1),\ldots, -1,1, \ldots, n-1,n\}.\)

Finally, we briefly recall the definitions of a topological Hopf
algebra and a deformation of a Hopf
algebra. See \cite{CP} for details.

\begin{definition}
A {\em topological Hopf algebra over $\C[[h]]$} is a complete
$\C[[h]]$-module $A_h$ equipped with $\C[[h]]$-linear maps \(\imath_h,
\mu_h, \epsilon_h, \Delta_h, S_h\) satisfying the Hopf algebra axioms,
but with algebraic tensor products replaced by the completions in the
$h$-adic topology. 
\end{definition}

\begin{definition}
A {\em deformation} of a Hopf algebra \((A, \imath, \mu, \epsilon, 
\Delta, S)\) over $\C$ is a topological Hopf algebra $A_h$ 
over $\C[[h]]$ such that 
\begin{enumerate}
\item $A_h$ is isomorphic to $A[[h]]$ as a $\C[[h]]$-module, 
\item \(\mu_h \equiv \mu \quad (\mathrm{mod} h), \quad \quad 
\Delta_h \equiv \Delta
\quad (\mathrm{mod} h)\).
\end{enumerate}
\end{definition}
In particular, the compatibility condition that $\Delta_h$ is an
algebra homomorphism from \(A_h\) to \(A_h \otimes A_h\) is expressed
by 
\begin{equation}\label{eq:DeltahCompatible}
\Delta_h(\mu_h(a_1 \otimes a_2)) = (\mu_h \otimes
\mu_h)\Delta^{13}_h(a_1)\Delta^{24}_h(a_2).
\end{equation}
Here, for \(\Delta\) any coalgebra structure, we define $\Delta^{13}$
and $\Delta^{24}$ as follows. If \(\Delta(a) = \sum a_i \otimes a'_i,
\Delta(b) = \sum b_j \otimes b'_j,\) then 
\[
\Delta^{13}(a)\Delta^{24}(b) := \sum a_i \otimes b_j \otimes a'_i
\otimes b'_j.
\]
We are interested in the first-order terms in $h$. In particular, the
first-order term in $h$ of the 
above equation~\eqref{eq:DeltahCompatible} for $\Delta_h = \sum_{k=0}^{\infty} \Delta_k h^k$ and
$\mu_h = \sum_{k=0}^{\infty} \mu_k h^k$ on $A[[h]]$
is 
\begin{equation}\label{eq:Hopfcompatibility}
\Delta_1(\mu_1(a_1\otimes a_2)) + \Delta_1(\mu_0(a_1 \otimes a_2)) =
(\mu_0 \otimes \mu_1 + \mu_1 \otimes
\mu_0)\Delta^{13}(a_1)\Delta^{24}(a_2) + \Delta_1(a_1)\Delta_0(a_2) +
\Delta_0(a_1)\Delta_1(a_2).
\end{equation}

\bibliographystyle{plain}
\bibliography{ref}

\begin{thebibliography}{10}

\bibitem{CP}
V.~Chari and A.~Pressley.
\newblock {\em A Guide to Quantum Groups}.
\newblock Cambridge University Press, 2nd edition, 1994.

\bibitem{Cur}
M.~Curtis.
\newblock {\em Matrix Groups}.
\newblock Springer-Verlag, 2 edition, 1984.

\bibitem{Dix74}
J.~Dixmier.
\newblock {\em Alg\`ebres Enveloppantes}.
\newblock Gauthier-Villars (Paris), 2 edition, 1974.

\bibitem{GelTse2}
I.~Gel$'$fand and M.~Tsetlin.
\newblock Finite-dimensional representations of the group of orthogonal
  matrices.
\newblock {\em Dokl. Akad. Nauk SSSR}, 17:1017--1020, 1950.

\bibitem{GelTse1}
I.~Gel$'$fand and M.~Tsetlin.
\newblock Finite-dimensional representations of the group of unimodular
  matrices.
\newblock {\em Dokl. Akad. Nauk SSSR}, 71:825--828, 1950.

\bibitem{GS82}
V.~Guillemin and S.~Sternberg.
\newblock Geometric quantization and multiplicities of group representations.
\newblock {\em Invent. Math.}, 67(3):515--538, 1982.

\bibitem{GS83}
V.~Guillemin and S.~Sternberg.
\newblock The {G}el{$'$}fand-{C}etlin system and quantization of the complex
  flag manifolds.
\newblock {\em Journal of Function Analysis}, 52:106--128, 1983.

\bibitem{Kos65}
B.~Kostant.
\newblock Orbits, symplectic structures, and representation theory.
\newblock {\em Proc. US-Japan Seminar in Differential Geometry}, 1965.

\bibitem{LMTW}
E.~Lerman, E.~Meinrenken, S.~Tolman, and C.~Woodward.
\newblock Non-abelian convexity by symplectic cuts.
\newblock {\em Topology}, 37(2):245--259, 1998.

\bibitem{M1}
A.~Molev.
\newblock Gelfand-{T}setlin basis for representations of {Y}angians.
\newblock {\em Lett. Math. Phys.}, 30(1):53--60, 1994.

\bibitem{M4}
A.~Molev.
\newblock A basis for representations of symplectic {L}ie algebras.
\newblock {\em Comm. Math. Phys.}, 201(3):591--618, 1999.

\bibitem{MNO}
A.~Molev, A.~Nazarov, and G.~Ol{$'$}shanskii.
\newblock Yangians and classical {L}ie algebras.
\newblock {\em Russian Math. Surveys}, 51(2):205--282, 1996.

\bibitem{MO}
A.~Molev and G.~Ol{$'$}shanskii.
\newblock Centralizer construction for twisted {Y}angians.
\newblock {\em Selecta Math. (N.S.)}, 6(3):269--317, 2000.

\bibitem{O92}
G.~Ol{$'$}shanskii.
\newblock Twisted {Y}angians and infinite-dimensional classical {L}ie algebras.
\newblock {\em Quantum Groups (Lecture Notes in Mathematics)}, 1510:103--120,
  1992.

\end{thebibliography}

\end{spacing}

\end{document}